\newif\ifdraft
\drafttrue
 
\documentclass[3p, times, review ]{elsarticle}

\usepackage[utf8]{inputenc}
\usepackage[english]{babel}
  \usepackage[]{algorithm2e}
\usepackage[utf8]{inputenc}
\usepackage[T1]{fontenc}
\usepackage{array}
\usepackage{amssymb}
\usepackage{amsfonts}
\usepackage{amsmath}
\usepackage{cases}
\usepackage{hyperref}

\usepackage{amsthm}
\usepackage{graphicx}
\usepackage{subfig} 
\usepackage{textcomp}
\usepackage{bm}
\usepackage{booktabs}
\usepackage{mathtools}

\makeatletter

\newcommand{\Rmnum}[1]{\expandafter\@slowromancap\romannumeral #1@}
\makeatother

\theoremstyle{definition}

\newcommand{\Dox}{D_\mathrm{ox}}
\newcommand{\Dliq}{D_\mathrm{liq}}
\newcommand{\blue}[1]{\textcolor{black}{#1}}

\newcommand{\unit}[2]{#1\,\mathrm{#2}}

\usepackage{color}

\ifdraft
\usepackage{color}
\usepackage[draft,disable]{todonotes}

\else
\newcommand{\todo}[1]{}

\fi

\newcommand{\RN}[1]{%
  \textup{\uppercase\expandafter{\romannumeral#1}}%
}
 
\journal{Accepted for Publication in Journal of Computational Electronics}
 
\usepackage[figuresright]{rotating}

\begin{document}

\begin{frontmatter}


\title{  Bayesian inversion for nanowire field-effect sensors}
\author[tuwien,Hannover]{Amirreza Khodadadian\corref{cor1}}
\ead{khodadadian@ifam.uni-hannover.de}
\author[tuwien]{Benjamin   Stadlbauer}
\ead{benjamin.stadlbauer@tuwien.ac.at}
\author[tuwien,ASU]{Clemens Heitzinger}
\ead{clemens.heitzinger@tuwien.ac.at}
\cortext[cor1]{Corresponding author}

\address[tuwien]{Institute for Analysis and Scientific Computing, 
  Vienna University of Technology (TU Wien),
  Wiedner Hauptstraße 8--10, 
  1040 Vienna, Austria}
  \address[Hannover]{
  	Institute of Applied Mathematics, Leibniz University Hannover, Welfengarten 1, 30167
  	Hanover, Germany
  }

  \address[ASU]{School of Mathematical and Statistical Sciences, Arizona State University, Tempe, AZ 85287, USA}


\begin{abstract}
  Nanowire field-effect sensors have  recently been developed for label-free detection of biomolecules. In this work, we introduce a computational technique based on Bayesian estimation to determine the physical parameters of the sensor and, more importantly, the properties of the analyte molecules. To that end, we first propose a PDE based model to simulate the device charge transport and electrochemical behavior. Then, the adaptive Metropolis algorithm with delayed rejection (DRAM)  is applied to estimate the posterior distribution of unknown parameters, namely molecule charge density, molecule density, doping concentration, and electron and hole mobilities.  We determine the device and molecules properties simultaneously, and we also calculate the molecule density as the only parameter after having determined the device parameters. This approach makes it possible not only to determine unknown parameters, but it also shows how well each parameter can be determined by yielding the probability density function (pdf).

\end{abstract}


\end{frontmatter}

\section{Introduction}\label{s:intro}
 
\begin{figure}[ht!]
    \centering
    \includegraphics[width=14cm,height=10cm]{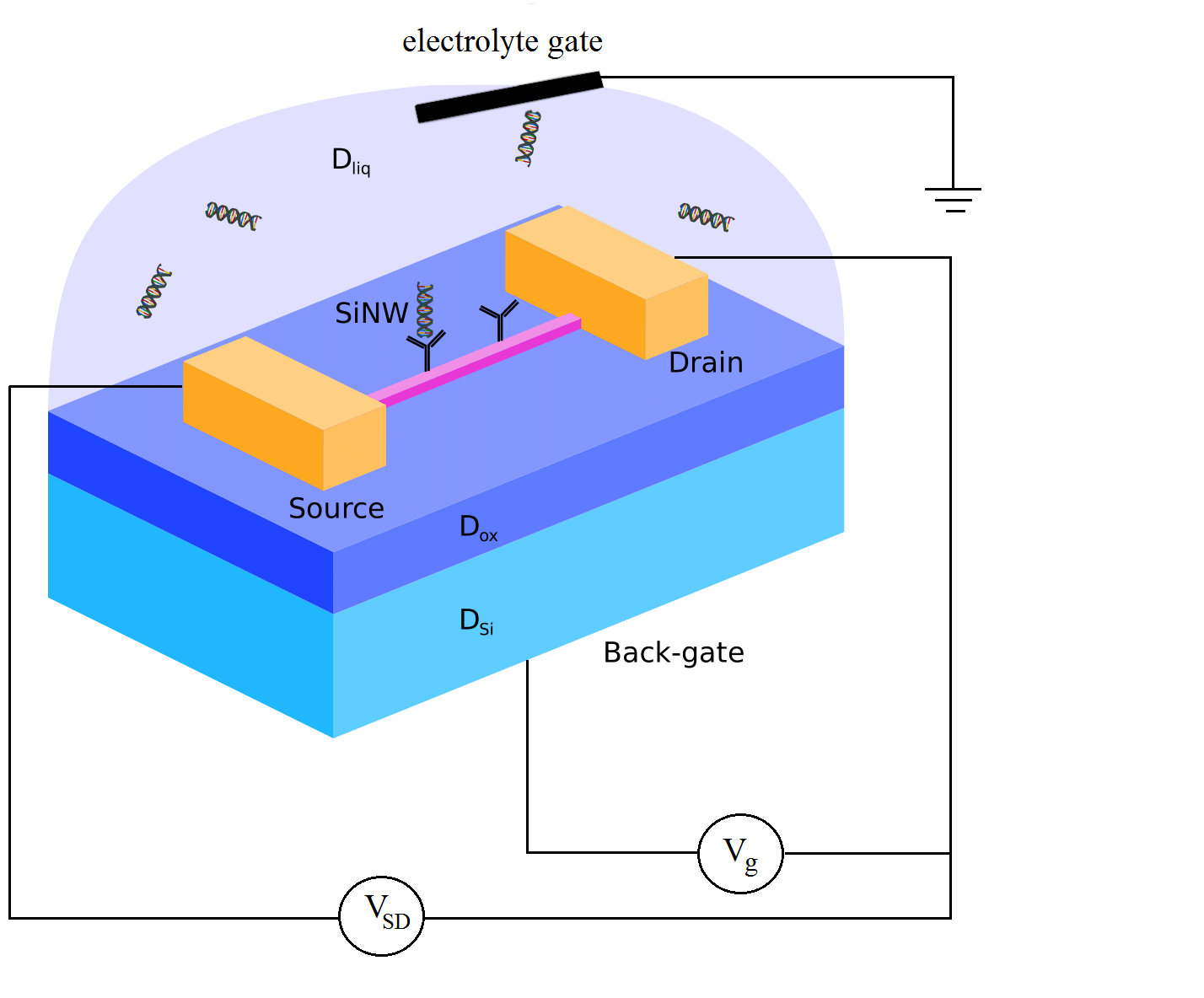}
    \caption{A schematic diagram of a SiNW sensor. The device consists of three  subdomains, i.e., silicon nanowire ($D_\text{Si}$), silicon dioxide insulator ($D_\text{ox}$) and electrolyte ($D_\text{liq}$). The nanowire is coated with a thin oxide layer. The target molecules and their binding to the receptors is illustrated as well.}
    \label{fig:schematic}
  \end{figure}

 Silicon nanowire sensors (SiNW) \cite{duan2013complementary,cui2001nanowire} (see Figure \ref{fig:schematic}) are promising devices used to detect 
 the presence or concentration of
  different biological species, such as cancer cells \cite{zheng2005multiplexed}, DNA and miRNA molecules \cite{he2015label,hahm2004direct}, and proteins \cite{wang2005label}. The sensors are being developed for the early detection of cardiovascular diseases \cite{chua2009label}, prostate cancer \cite{baumgartner2013predictive}, breast cancer \cite{shashaani2016silicon}, gastric cancer \cite{lee2010measurements}, flu \cite{shen2012rapid}, and uric acid in human blood \cite{guan2014highly}.  
  In the sensors, 
  the target molecules such as
  biomarkers bind selectively to the recognition elements, e.g., antibodies or aptamers.
  The semiconductor transducer
  converts the potential change due to the analyte molecules into a measurable
  electrical signal, i.e., a current or voltage change \cite{mu2015silicon}.  The biosensors are interesting candidates for  biomarker detection since the sensors are reliable, label-free, inexpensive,  highly sensitive, and have short operation time \cite{khodadadian2017optimal, mirsian2019new}. 
 
The binding of target molecules to the receptors changes the charge concentration, which modulates the transducer. The drift-diffusion-Poisson-Boltzmann (DDPB) system is a comprehensive set of equations to model the device. The Poisson-Boltzmann equation enables us to calculate the charge concentrations around biomolecules and the effect of the charged target
molecules on the nanowire, while the drift-diffusion equations model the charge transport of the carriers through the nanowire.  In the biosensors, there are many sources of noises, e.g., random movement and binding of the target molecules \cite{khodadadian2017optimal,khodadadian2016basis}. There are also random dopant fluctuations (RDF) \cite{taghizadeh2017optimal}. The stochastic version of the equation system can be used to investigate these effects. Also, existence and uniqueness theorems for the stochastic equations are given in \cite{taghizadeh2017optimal}.  

For solving the stochastic DDPB system, efficient numerical methods are finite elements (for space discretization) and the Monte Carlo  technique (for the stochastic dimensions). In \cite{taghizadeh2017optimal,khodadadian2018three}, the authors developed a multilevel Monte Carlo  finite element method (MLMC-FEM) to solve the equation with little computational effort. A further complexity reduction has been been achieved by replacing the random points by quasi-Monte Carlo points \cite{khodadadian2018optimal}. Also, by using three-dimensional simulations \cite{khodadadian2018three,baumgartner2013predictive} we can perform realistic simulations.

In order to obtain good agreement with experiments, a robust estimation of unknown model parameters, e.g., diffusion coefficients, the charge and density of molecules, doping concentration, etc., is essential. In classical inversion methods, one strives to minimize the distance between measurement and simulation to estimate the unknown parameters.
Minimizing the difference between measurement and simulation always yields a value, but its sensitivity must be assessed separately. This is not a trivial manner especially in the case of ill-posed nonlinear problems. The great advantage of the Bayesian approach followed here is that it yields the probability density of the unknown parameters. 


 In Bayesian inversion \cite{dashti2017bayesian} the solution of the inverse problem is the posterior density giving the distribution of the unknown parameter values based on the sampled observations \cite{smith2013uncertainty}.
 Markov-chain Monte Carlo (MCMC) is a popular method to calculate the distribution. In this method, a Markov chain is constructed whose stationary distribution is the sought posterior distribution in Bayes theorem.
 

The Metropolis-Hastings (MH) algorithm \cite{smith1993bayesian} is one of the most common techniques among the MCMC methods since it is simple and sufficiently powerful for many problems (specifically when the parameters are not strongly correlated). In order to estimate the posterior distribution, in each iteration, we propose a new candidate  parameter value based on the current sample value according to a proposal distribution. Then, we calculate the acceptance ratio and decide whether the candidate value is accepted or rejected.  The acceptance ratio points out how probable the new candidate value is with respect to the current sample. 

 The MH algorithm has some drawbacks, e.g., the proposal covariance must be manually tuned and has high autocorrelation   \cite{andrieu2003introduction}. To overcome these deficiencies, instead of using a fixed proposal distribution in each iteration, we update the distribution according to the available samples (adaptive Metropolis). This approach is useful since the posterior distribution is not sensitive to the proposal distribution. 
 The adaptive method can be modified additionally by combining it with delayed rejection yielding the DRAM algorithm. In this algorithm, an alternative for the rejected candidate is proposed and the probability of this conditional acceptance is corrected \cite{mira2001metropolis}. Upon rejection in the
 MH algorithm, instead of  retaining the current position, a second-stage move is
 proposed \cite{haario2006dram}. The method is noticeably advantageous when  sampling from high-dimensional conditional distributions \cite{zuev2011modified}.

%

In the sensor design process, we use the PDE-based model to simulate the device characteristic. However, a very good estimation of the model parameters enables us to have a more efficient simulation and predict the sensor electrical behavior in different situations (e.g., subthreshold and linear regimes). In this work, the main aim is to propose an efficient computational method to estimate physical parameters of the sensors that cannot be measured directly or only with great experimental efforts \cite{punzet2012determination}. For instance, regarding the target molecules, a reliable estimation of the surface charge (due to their binding to the receptors) or their reaction with the probe molecules cannot be achieved easily. Similarly, a reasonable estimation of the number of target molecules bound to the receptors cannot be achieved experimentally. However, in order to have an exact simulation, these parameters are crucial.
To that end, we use the DRAM algorithm to calculate the posterior distribution of various unknown parameters. The extracted information will help us to improve the simulation quality since as a biosensor and a transistor the necessary parameters are determined efficiently.

  The outline of the paper is as follows. In Section \ref{section2}, we introduce the stochastic drift-diffusion-Poisson-Boltzmann model and explain how it can be used to model all charge interactions in the device. In Section \ref{section3}, we present the MH algorithm and its modification by covariance adaption and delayed rejection. In Section \ref{section4}, we use  the PDE model and the DRAM technique to estimate the unknown parameters of a nanowire sensor. In order to do this, we first validate the transport model by comparison with experimental data and then 
 use these experimental data in the Bayesian inversion. Finally, the conclusions are drawn in Section \ref{conclusions}.

\section{The macroscopic model equation}
\label{section2}


In this section, we review a complete mathematical model to understand the physical behavior of nanowire sensors where the details are explained in the authors previous papers, e.g., \cite{baumgartner2013predictive,khodadadian2018optimal,khodadadian2017optimal}.  
It
consists of the Poisson-Boltzmann equation to model the electrolyte and the drift-diffusion-Poisson
system to model the charge transport in the semiconducting part.
 
Here we assume $D\subseteq \mathbb{R}^3$ is the sensor domain (our computational geometry) which is partitioned into three regions
with their physical characteristics.
The subdomain $D_{\mathrm{Si}}$ (silicon nanowire)
is the transducer insulated by $\text{SiO}_2$ ($D_{\mathrm{ox}}$), the second subdomain. In $D_{\text{liq}}$, the electrolyte contains cations and anions; therefore, the Poisson-Boltzmann equation holds.  

 In order to describe the sensor electrostatic interactions, we use the stochastic Poisson-Boltzmann equation 
  \label{modeleqn1}
		\begin{alignat}{3}
		-\nabla\cdot(A(x)\nabla  V(x,\omega))&= 
		\begin{cases}\label{model} 
		q(C_\mathrm{dop}(x,\omega)+p(x,\omega)-n(x,\omega))&  \text{in} \; D_{\text{Si}}, \\     
		0 & \text{in} \; D_{\mathrm{ox}}, \\  
		-2\varphi(x,\omega) \sinh (\beta(V(x,\omega)- \Phi_F))+ \rho(x,\omega) & \text{in} \;D_{\text{liq}},
		\end{cases}
		  \end{alignat}
		 where $A$ is the  dielectric (permittivity) function with the relative
		 permittivities of the materials  assumed to be constant and equal to $A_{\text{Si}}=11.7$,
		 $A_{\text{ox}}=3.9$, and $A_{\text{liq}}=78.4$.  In this equation, $x\in D$ and $\omega$ belongs to $\Omega$  which is the probability space. Moreover,
		  $V$ is the electrical potential, $q$ is the elementary charge, $C_\text{dop}$ is the doping concentration, $\rho$ considers the molecules charge, $\varphi$ is the ionic concentration (holds for a  symmetric electrolyte of monovalent ions), $\beta:=q/K_BT$,   and  $\Phi_F$ is the Fermi level. The concentration of electrons and holes are given by a Boltzmann distribution as
\begin{align}
p(x,\omega)&=n_i\exp\left(-\frac{qV(x,\omega)-\Phi_F}{K_BT}\right),\\
n(x,\omega)&=n_i\exp\left(\frac{qV(x,\omega)-\Phi_F}{K_BT}\right),
\end{align}
where $K_B$ is the Boltzmann constant and $T$ is the temperature.

The interface conditions in the electrostatic potential arise from homogenization \cite{heitzinger2010multiscale}. The interface conditions are
  \begin{align}
  V(0+,y,\boldsymbol{\omega}) - V(0-,y,\boldsymbol{\omega})
  & = \blue{\frac{D}{A^+}}
  && \text{on } \Gamma,\\
  A(0+)\partial_x V(0+,y,\boldsymbol{\omega})-A(0-)\partial_x V(0-,y,\boldsymbol{\omega})
  & = \blue{-C}
  && \text{on } \Gamma,
  \end{align}
  where the interface is $\Gamma=\Dliq\cap\Dox$.  \blue{Here, $0+$ indicates
  the limit at the interface on the outside of the sensor, while
  $0-$ is the limit on the inside. The two interface conditions
  contain the cumulative effect of a rapidly oscillating charge
  concentration in the surface or boundary layer at the surface of the
  sensor.  The constant $A^+$ is the permittivity outside the sensor,
  i.e., the permittivity of $\Dliq$.  The constant~$C$ is the
  macroscopic surface-charge density of the boundary layer, and the
  constant~$D$ is its macroscopic dipole-moment density as defined in
  \cite{heitzinger2010multiscale}.}
  		 	
  For this model, the boundary conditions  are Dirichlet boundary conditions ($\partial D_D$) and Neumann boundary condition ($\partial D_N$). A voltage across the simulation domain in the
  	vertical direction can be applied   by an electrode in the
  	liquid (solution voltage) and by a back-gate contact at the bottom
  	of the structure (back-gate voltage). At the source and drain contacts, the Dirichlet boundary conditions are respectively $V_S$ (source voltage) and $V_D$ (drain voltage). Additionally, zero Neumann boundary conditions hold on else everywhere (the
  	Neumann part $\partial D_N$ of the boundary).
  
Next, we consider the drift-diffusion-Poisson equations to model the charge transport through the semiconducting nanowire. In the transducer $D_{\mathrm{Si}}$, the stochastic drift-diffusion-Poisson system \cite{taghizadeh2017optimal}
\begin{subequations}
\begin{align}
-\nabla\cdot(A(x)\nabla  V(x,\omega))&=q(C_\mathrm{dop}(x,\omega)+p(x,\omega)-n(x,\omega)),\ 
\\ 
\nabla \cdot J_n(x,\omega)&=qR(n(x,\omega),p(x,\omega)),\ 
\\ 
\nabla\cdot J_p(x,\omega)&=-qR(n(x,\omega),p(x,\omega)),\ 
\\ 
J_n(x,\omega)&=q(D_n\nabla n(x,\omega)-\mu_nn(x,\omega)\nabla V(x,\omega)),\ 
\\ 
J_p(x,\omega)&=q(-D_p\nabla p(x,\omega)-\mu_pp(x,\omega)\nabla V(x,\omega))
\end{align}
\label{ddp}
\end{subequations}
\hspace{-.1cm}holds. Here,  $J_n(x,\omega)$ and $J_p(x,\omega)$ indicate the current
densities of the carriers, $D_n$ and $D_p$ are the diffusion coefficients, $\mu_n$ and $\mu_p$ are the mobilities, and 
$R(n(x,\omega),p(x,\omega))$ is the recombination rate. We use the Shockley-Read-Hall (SRH) recombination rate, which is defined as
\begin{align}\label{equation2}
R(n(x,\omega),p(x,\omega)):=\frac{n(x,\omega)p(x,\omega)-n_i^2}{\tau_p(n(x,\omega)+n_i)+\tau_n(p(x,\omega)+n_i)},
\end{align}
where $n_i:=1.5 \times 10^{10}\mathrm{cm^{-3}}$ is the intrinsic charge density and $\tau_n$ and $\tau_p$ are the lifetimes of the free
carriers. 
Finally, the total current density is $J_n+J_p$ and the total electrical current
\begin{align}\label{current}
\mathcal{I}:=\int \left(J_n+J_p\right) \,\text{d}x
\end{align}
 is obtained by calculating this integral over a cross-section of the transducer.

 In order to include the biological noise, we consider an association/dissociation process \cite{tulzer2014fluctuations}. The association and dissociation processes of target molecules at the surface can
 be described by the reaction process
 \begin{subequations} 
 	\label{reaction}
 	\begin{align}
 	\boldsymbol{\mathrm{T}}+ \boldsymbol{\mathrm{P}}&~{\longrightarrow}~\boldsymbol{\mathrm{PT}}, \\
 	\boldsymbol{\mathrm{PT}}&~{\longrightarrow}~\boldsymbol{\mathrm{P}}+\boldsymbol{\mathrm{T}},
 	\end{align}
 	\label{PT}
 \end{subequations}
\hspace{-.15cm} where (\ref{PT}) describes association and dissociation of the probe-target complex at the surface. In other words, the binding of target molecules $\boldsymbol{\mathrm{T}}$   (target-molecule
 concentration) to probe molecules $\boldsymbol{\mathrm{P}}$ (probe-molecule
 concentration), thus forming the probe-target complex $\boldsymbol{\mathrm{PT}}$ (probe-target
 concentration). The reaction equations provides sufficient information about the number of bound molecules to the receptors ($\boldsymbol{\mathrm{PT}}$-complex at the surface).

\section{Metropolis-Hastings algorithm}
\label{section3}
In this section, we briefly introduce the Bayesian inversion approach and explain how this technique can be implemented to estimate the unknown parameters. First we use the statistical model
\begin{align}
\mathcal{M}_i=\mathcal{I}_i(A)+\varepsilon_i,\qquad i=1,\ldots,n,
\end{align} 
where $\mathcal{M}_i$, $\mathcal{I}_i,$ and $\varepsilon_i$ are random
variables representing the measurement, the estimated current by the
model (here the drift-diffusion-Poisson system (\ref{current})), and
the measurement error, respectively. The measurement error is a
realization of $\mathcal{N}(0,\sigma^2 I)$, where \blue{$\mathcal{N}$
  indicates the normal distribution, $I$ is the identity matrix,} and $\sigma^2$, its variance, is
a fidelity parameter that corresponds to the measurement error.
 For a given value $\alpha$ of the parameter $A$ and the corresponding observation $\beta$ (the measurement)
  we assume that
  $\pi$ is a Lebesgue density. Then, we define the conditional density
  \begin{align*}
  \pi(\alpha|\,\beta)=\frac{\pi(\alpha,\beta)}{\pi(\beta)}
  \end{align*}      
  such that
  \begin{eqnarray}
\pi(\beta) =\int_{\mathbb{R}^{d_1}}\pi(\beta,\alpha)~\pi_0(\alpha)\,\text{d}\alpha,
  \end{eqnarray}
  where $\pi_0(\alpha)$ is the $\textit{prior}$ probability density \cite{stuart2010inverse}. Using the measured value $\beta$, Bayes Theorem yields the \textit{posteriori} distribution as
  \begin{align}\label{post}
\pi(\alpha|\,\beta)=\frac{\pi(\beta|\,\alpha)\pi_0(\alpha)}{\pi(\beta)}=\frac{\pi(\beta|\,\alpha)\pi_0(\alpha)}{\int_{\mathbb{R}^{d_1}}\pi(\beta,\alpha)~\pi_0(\alpha)\,\text{d}\alpha }.
  \end{align} 
 In \eqref{post} the denominator is a normalization constant and its explicit calculation is computationally expensive. Therefore, we sample the posterior distribution without the knowledge of the normalization constant yielding 
\begin{align*}
\pi(\alpha|\,\beta)\propto\pi(\beta|\,\alpha)\pi_0(\alpha).
\end{align*} 
Now we assume that the errors are \blue{independent and identically
  distributed} (iid) and that $\varepsilon_i\sim N(0,\sigma^2)$. The \blue{likelihood function}  $\blue{\pi(\beta|\,\alpha)}$ is therefore 
\begin{align}\label{likelihood}
\pi(\beta|\,\alpha):=\frac{1}{\left(2\pi\sigma^2\right)^{n/2}}\exp(-\zeta(\alpha)/2\sigma^2),
\end{align}
where
\begin{align}
\label{errr}
\zeta(\alpha)=\sum_{i=1}^{n}|\,\beta_i-\mathcal{I}_i(\alpha)\,|^2,
\end{align}
indicates the simulation error with respect to the parameter
$\alpha$. \blue{Here \eqref{errr} indicates the difference between the
  measurement $\beta_i$ and the solution $\mathcal{I}_i$ of the forward model (i.e.,
  the DDP system) for the proposed candidate $\alpha$ for the
  different cases $i \in \{1,\ldots,n\}$. Therefore, the likelihood
  function $\pi$ expresses the plausibilities of different parameter
  values $\alpha$ given the observations $\beta$.}

 When we have enough information about the posterior distribution in
the most straightforward situation, we can directly sample from it. However, in most   cases, we do not have sufficient knowledge about the distribution or it is not possible to sample from it due to high-dimensionality or complexity. To overcome this problem, the MH algorithm can be  used. A summary of the algorithm is given in Algorithm 1.

\begin{algorithm}[ht!]
   \label{algorithm}
   \textbf{Initialization ($\ell=0$)}: Generate the initial parameter $\alpha^0\sim\pi(\alpha^0|\,\beta)$.\\

   \For{$\ell<N$}{

     \quad 1. $\ell=\ell+1$.\\
 
    \quad 2. Propose the new candidate $\alpha^{*}\sim\phi(\alpha^\ell|\,\alpha^{\ell-1})$.    \\

    \quad 3. Compute the proposal correction parameter~$\upsilon(\alpha^*|\,\alpha^{\ell-1})=\cfrac{\phi(\alpha^{\ell-1}|\,\alpha^*)}{\phi(\alpha^*|\,\alpha^{\ell-1}))}$.
    \\

   \quad  4. Calculate the acceptance/rejection probability $\lambda (\alpha^*|\,\alpha^{\ell-1})=\min\left(1, \cfrac{\pi(\alpha^*|\,\beta)}{\pi(\alpha^{\ell-1}|\,\beta)}\,\upsilon\right)$.\\
 
       \quad 5. Draw a random number $\mathcal{R}\sim \text{Uniform}\,(0,1)$.\\
 
   \quad 6.$~\textbf{if}~$ $\mathcal{R}<\lambda~$ \textbf{then}\\ 
 
   \qquad\quad accept the candidate $\alpha^*$ and set $\alpha^\ell=\alpha^*$\\
 
   \qquad \textbf{else}\\
 
         \qquad\quad reject the candidate $\alpha^*$ and set $\alpha^{\ell}=\alpha^{\ell-1}$\\
  
   \qquad \textbf{end if}

   }
   \caption{The MH algorithm.}
 \end{algorithm}
   In the algorithm, the proposal distribution \blue{$\phi$} is fixed; therefore,  the rejection rate can be very high. Using an updated covariance matrix for the proposal distribution (applying the learned information about the posterior) allows us to increase the acceptance ratios since they accelerate the rate at which information regarding the posterior is incorporated \cite{smith2013uncertainty}. Also, to enhance  efficiency, the adaptive algorithm is combined with a delayed rejection technique. To that end, using the information of the rejected proposal, a new candidate is proposed and is rejected or accepted based on a suitably computed probability \cite{haario2006dram}. A summary of the DRAM method is given in Algorithm 2.
   
   \begin{algorithm}[ht!]
      \label{algorithm1}
      \textbf{Initialization ($\ell=0$)}: Generate the initial parameter $\alpha^0\sim\pi(\alpha^0|\,\beta)$.\\

      \For{$\ell<N$}{

        \quad 1. $\ell=\ell+1$.\\
       
       \quad 2. Propose the new candidate $\alpha^{*}=\alpha^{\ell-1}+\mathcal{R}_\ell\mathcal{Z}_\ell$  where $\mathcal{R}_\ell$ is the Cholesky decomposition of $\mathcal{V}_\ell$.  \\

      \quad  3. Calculate the acceptance/rejection probability $\lambda_1 (\alpha^*|\,\alpha^{\ell-1})=\min\left(1, \cfrac{\pi(\alpha^*|\,\beta)~\phi(\alpha^{\ell-1}|~\alpha^*)}{\pi(\alpha^{\ell-1}|\,\beta)~\phi(\alpha^*|\,\alpha^{\ell-1}))}\right)$.\\
 
          \quad 4. Draw a random number $\mathcal{R}\sim \text{Uniform}\,(0,1)$.\\
     
      \quad 5.$~\textbf{if}~$ $\mathcal{R}<\lambda_1~$ \textbf{then}\\ 
 
      \qquad\quad accept the candidate $\alpha^*$ and set $\alpha^\ell=\alpha^*$\\
 
      \qquad \textbf{else}\\
 
            \qquad\quad1.  Calculate the alternative candidate
            \quad$\alpha^{**}=\alpha^{\ell-1}+\sigma^2\mathcal{R}_\ell\mathcal{Z}_\ell$.\\
            \qquad\quad 2.   Calculate the acceptance/rejection probability\\ $\qquad\qquad\lambda_2 (\alpha^{**}|\,\alpha^{\ell-1},\alpha^*)=\min\left(1, \cfrac{\pi(\alpha^{**}|\,\beta)~\phi(\alpha^*|~\alpha^{**})\left(1-\lambda_1(\alpha^*|\alpha^{**})\right)}{\pi(\alpha^{\ell-1}|\,\beta)~\phi(\alpha^*|~\alpha^{\ell-1})\left(1-\lambda_1(\alpha^*|\alpha^{\ell-1})\right)}\right)$.\\
 
           \qquad\quad 3.  \textbf{if} ~$\mathcal{R}<\lambda_2$~~\textbf{then}\\
         \qquad\qquad\quad~~accept the candidate $\alpha^{**}$ and set $\alpha^{\ell}=\alpha^{**}$\\
         \qquad\quad   \textbf{else}\\
  \qquad\qquad\quad~~  reject the candidate $\alpha^{**}$ and set $\alpha^{\ell}=\alpha^{\ell-1}$
           \qquad\quad\\
      \qquad \textbf{end if}\\
    \quad 6. Update the covariance matrix as $\mathcal{V}_\ell=\text{Cov}(\alpha^0,\alpha^1,\ldots,\alpha^{\ell})$.\\
    \quad 7. Update $\mathcal{R}_\ell$.
      }
      \caption{The DRAM algorithm.}
    \end{algorithm}

  In Algorithm 2, $\mathcal{Z}_\ell\sim \text{Uniform}~(0,I_N)$ where $I_N$ is the $N$-dimensional identity matrix. In order to use a narrower proposal function (compared to the first proposal) we employ $\sigma<1$. The covariance function is calculated as
\begin{align}
\text{Cov}(\alpha^0,\alpha^1,\ldots,\alpha^{\ell})=\frac{1}{\ell}\left(\sum_{i=0}^{\ell}
\alpha^i\left(\alpha^i\right)^T-(\ell+1)~\hat{\alpha}^\ell\left(\hat{\alpha}^\ell\right)^T\right),
\end{align}  
where $\hat{\alpha}^\ell=\frac{1}{\ell+1}\sum_{i=0}^{\ell}\alpha^i$. We refer the interested reader to \cite{haario2006dram} for more details.

The main aim of using Bayesian inversion in this work is matching the electrical current $\mathcal{I}$ obtained by the drift-diffusion model \eqref{current} by the experimental measurements at different gate voltages. To that end, we provide a list of desirable unknown parameters (with their relative equations) and explain why they should be determined precisely.  
\begin{itemize}
	\item Surface charge density ($\rho$), in the stochastic Poisson-Boltzmann equation describes the charge due to binding the target molecules to the receptors at the sensor surface, i.e., $-\nabla\cdot(A(x)\nabla  V(x,\omega))=\rho(x,\omega)$. A reliable estimation due to the unknown area of probe-target molecule (therefore the surface charge estimation) cannot be easily estimated.
	
	\item Doping concentration ($C_\text{dop}$) in the stochastic Poisson-Boltzmann equation and drift-diffusion equations (in $D_{\text{Si}}$) is another influential parameter. Generally, it is an average number of dopants, and the precise concentration is not extracted usually. The exact doping density makes the simulation more reliable.   
	
	\item Electron mobility ($\mu_n$) and hole density ($\mu_p$) used in the stochastic drift-diffusion equation. As we already know, the sensor acts as a transistor; therefore, its electrical behavior in the subthreshold and linear regime is crucial. 
	Einstein relations, i.e., $D_n=U_T \mu_n$ and $D_p=U_T \mu_p$ ($U_T$ is the thermal voltage) gives us good information about the diffusion coefficients (used in the subthreshold conduction).
	 
	 \item $\textbf{PT}$-density presented in the reaction equations (\ref{reaction}) which give us the density of the probe-target complex at the sensor surface. Considering the surface area, we can predict how many target molecules absorb precisely to the receptors.
	  
\end{itemize}
 
 \blue{It is noted that there is no correlation between $\rho$,
   $C_\text{dop}$, and the $\textbf{PT}$-density; however, a higher
   doping concentration reduces the electron and hole
   mobility.  Parameter estimation using the DRAM algorithm and a set
   of backgate voltages $\mathcal{V}_\mathcal{G}$ is summarized in
   Algorithm \ref{algorithm3}.}

   \begin{algorithm}[ht!]
 	\textbf{Prior distribution}: Generate the initial parameter
 	 $\alpha^0\sim\pi(\alpha^0|\,\beta)$.\\  
 	 	
 	\For{$\ell=1:N$}{		
 		\vspace{0.1cm}
  	\quad 1. Propose  the $\ell$-th candidate $\alpha^*=(\rho^*,~ C_\text{dop}^*,~\mu_n^*,~\mu_p^*)$ according to the proposal distribution.    \\\vspace{0.1cm}
 		\quad	2. \textbf{for all} $V_g\in \mathcal{V}_\mathcal{G}$\\
 		\qquad\quad I.  Use Scharfetter-Gummel iteration \cite{taghizadeh2017optimal} to solve the forward model and estimate $J_n$ and $J_p$.\\\vspace{0.1cm}
 	 	\qquad\quad II. Take integral over a cross-section of semconducting part to obtain the electrical current \eqref{current}. \\\vspace{0.1cm}
 	 	\quad ~~~\textbf{end for}\\\vspace{0.1cm}
  	\quad 3. Calculate the likelihood function \eqref{likelihood} considering all $V_g\in \mathcal{V}_\mathcal{G}$  with respect to the measurements.\\\vspace{0.1cm}
  	\quad 4. Calculate the acceptance/rejection probability $\lambda_1$.\\\vspace{0.1cm}
 	 	\quad 5. \textbf{if} $\alpha^*$ is accepted \textbf{then}\\\vspace{0.1cm}
  	\qquad \quad set $\alpha^\ell =\alpha^*$\\\vspace{0.1cm}
  	\quad~~~ \textbf{else}\\\vspace{0.1cm}
 	 	\qquad \quad Calculate the alternative candidate $\alpha^{**}$ and its acceptance/rejection probability $\lambda_2$\\\vspace{0.1cm}
 	 	\qquad \quad Check that it is accepted/rejected and update the covariance (see DRAM Algorithm).	 \\\vspace{0.1cm}
  	\quad~~~  \textbf{end if}\\\vspace{0.1cm}
 	}
	\textbf{Posterior distribution:} Estimate the distribution according to all proposed candidates $\alpha$ and extract the information related to unknown parameters.\\\vspace{0.5cm}
 	\caption{The DRAM algorithm for silicon nanowire sensors using DDP system.} 
 	\label{algorithm3}	
 \end{algorithm}

\section{Parameter estimation and model verification}
\label{section4}
In this section, we use the DDPB system to model the electrical behavior of the sensor whose  3D schematic diagram is shown in Figure \ref{fig:schematic}. Regarding the geometry, the nanowire length is $\unit{1000}{nm}$, and its width and thickness are $\unit{50}{nm}$ and $\unit{100}{nm}$, respectively. The insulator (silicon dioxide) thickness is $\unit{8}{nm}$, the distance between the nanowire and the boundary is $\unit{50}{nm}$. In the simulations the thermal voltage $(1/\beta)$ is $\unit{0.021}{V}$ and the source-to-drain voltage is $\unit{0.2}{V}$.
The experimental data for the mentioned device are taken from \cite{baumgartner2013predictive} for the nanowire field-effect PSA sensor. Furthermore, we assume that the measurement has 1\% error.

\begin{figure}[ht!]
   \centering
   \includegraphics[width=10cm,height=7cm]{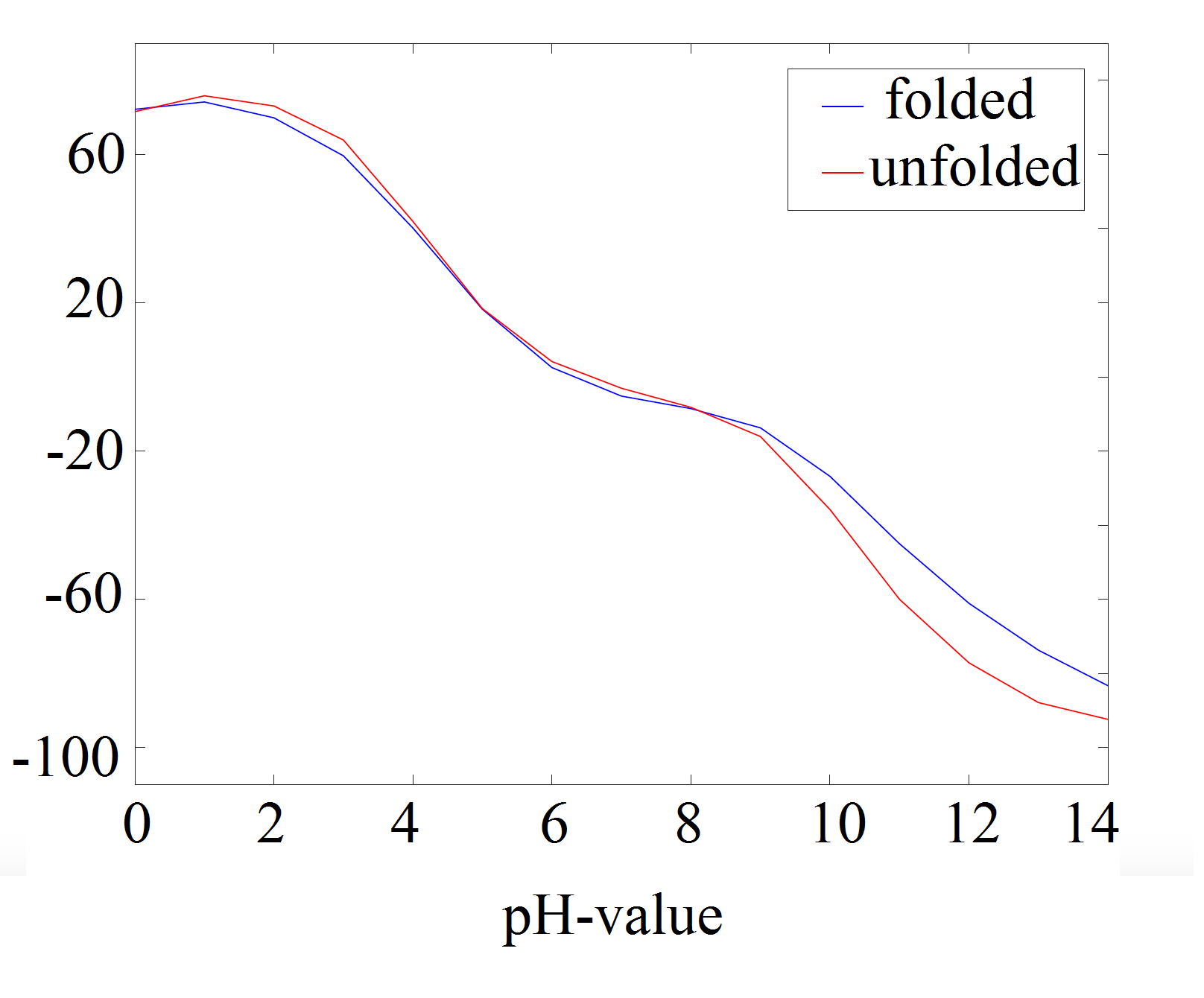}
   \caption{The net charge of 2ZCH protein for different pH values.}
   \label{fig:pH}
 \end{figure}

Prostate-specific antigen (PSA) is an important biomarker widely used to diagnose prostate cancer. Here we develop a sensor to detect the protein 2ZCH (\url{https://www.rcsb.org/structure/2ZCH}). 
The charge of PSA is a function of pH value, where the PROPKA algorithm \cite{li2005very} is applied to estimate the net charge shown in Figure \ref{fig:pH}. We perform the simulations with measurements performed at a pH value of 9; we note a total PSA charge of $\unit{-15}{q}$ at this pH value. 


 In P-type
semiconductors, applying a positive gate voltage depletes
carriers and reduces the conductance, while applying a
negative gate voltage gives rise to an accumulation of carriers 
and increases device conductivity. In field-effect biosensors, the PSA target molecules carry negative charges (see Figure \ref{fig:pH}) which act as a negative gate voltage.
 Since we use a P-type (boron-doped)
semiconductor as the transducer, the accumulation of holes increases the conductance
as well.

The doping concentration  is another important physical parameter in semiconductor devices.  In nanowire sensors, on the one hand, a higher concentration increases device conductivity, while, on
the other hand, a higher concentration decreases the sensor sensitivity. In other words, when the doping concentration is high, the nanowire is mostly affected by the dopant atoms and the effect of the charged molecules on the sensor response decreases.
Therefore, the optimal doping concentration in the device design process is essential.

\begin{figure}[ht!]
    \centering
    \includegraphics[width=10cm,height=6cm]{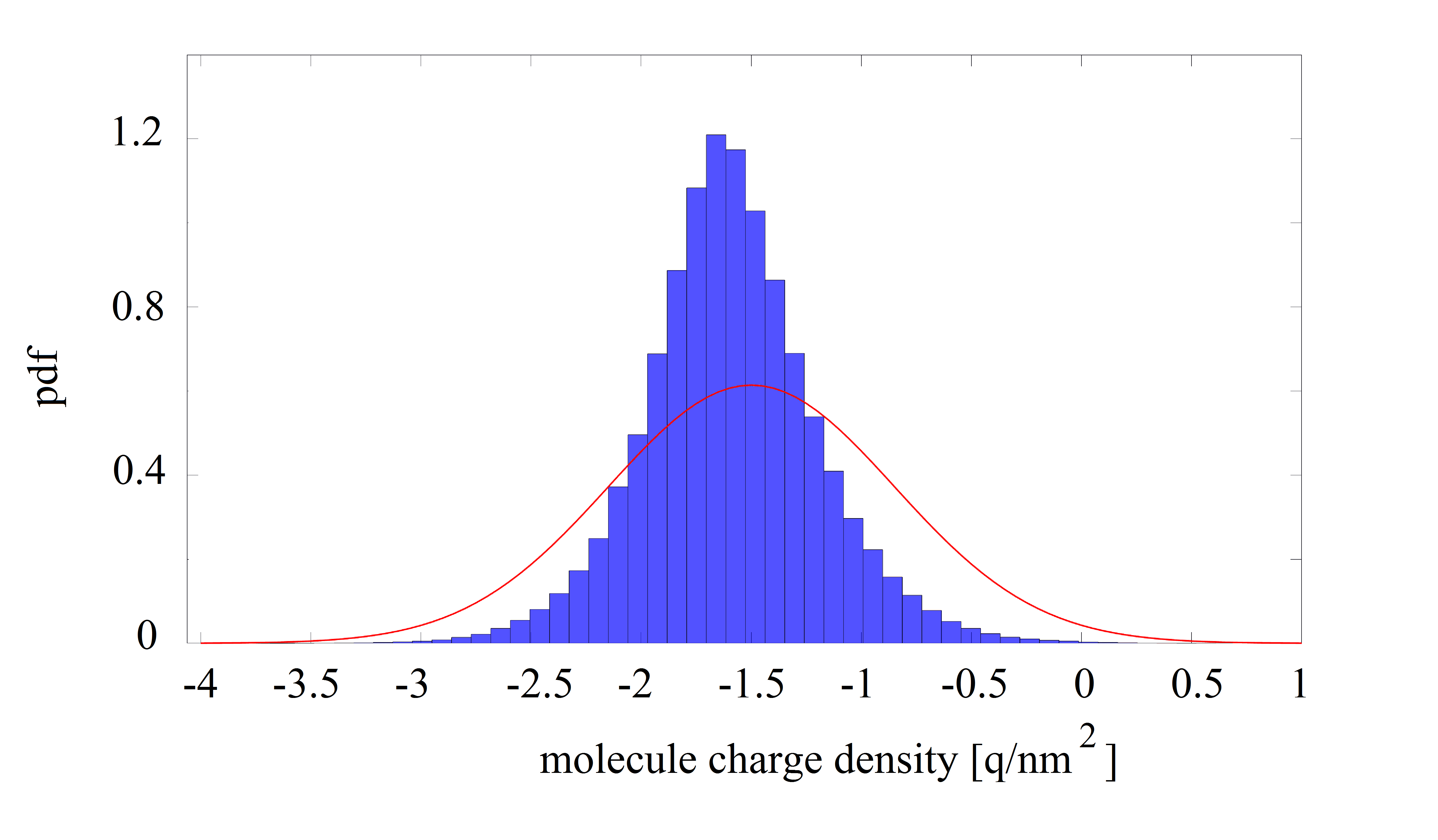}
    \caption{The probability density function (pdf) of posterior distribution (histogram) and prior distribution (red line) of the surface charge density using the DRAM algorithm.}
    \label{fig:surface1d}
  \end{figure}

From now on, we use the DRAM algorithm to estimate the important unknown parameters where in all cases $N=300\,000$ number of samples are used.
 The first study is the molecule charge density in the sense that the prior knowledge is $\rho=\unit{-1.5}{q/nm^2}$. Figure \ref{fig:surface1d} shows the posterior distribution and prior distribution  of the molecule charge density where the acceptance rate of 67.7\% is achieved. 
 Due to the obtained results by the PROPKA algorithm and using a P-type semiconductor we employ the (Gaussian) proposal distribution  between $\rho=\unit{1}{q/nm^2}$ and $\rho=\unit{-4}{q/nm^2}$.
 The results point out that most of the accepted proposals are around the prior knowledge (its expected value is $\rho=\unit{-1.71}{q/nm^2}$). The posterior distribution indicates that the probability of positive charges is negligible, which agrees very well with the transducer structure (P-type nanowire). Finally, the narrower shape of posterior distribution (compared to the proposal) indicates the Bayesian inversion efficiency.

\begin{figure}[ht!]
   \centering
   \subfloat{\includegraphics[width=8cm,height=5cm]{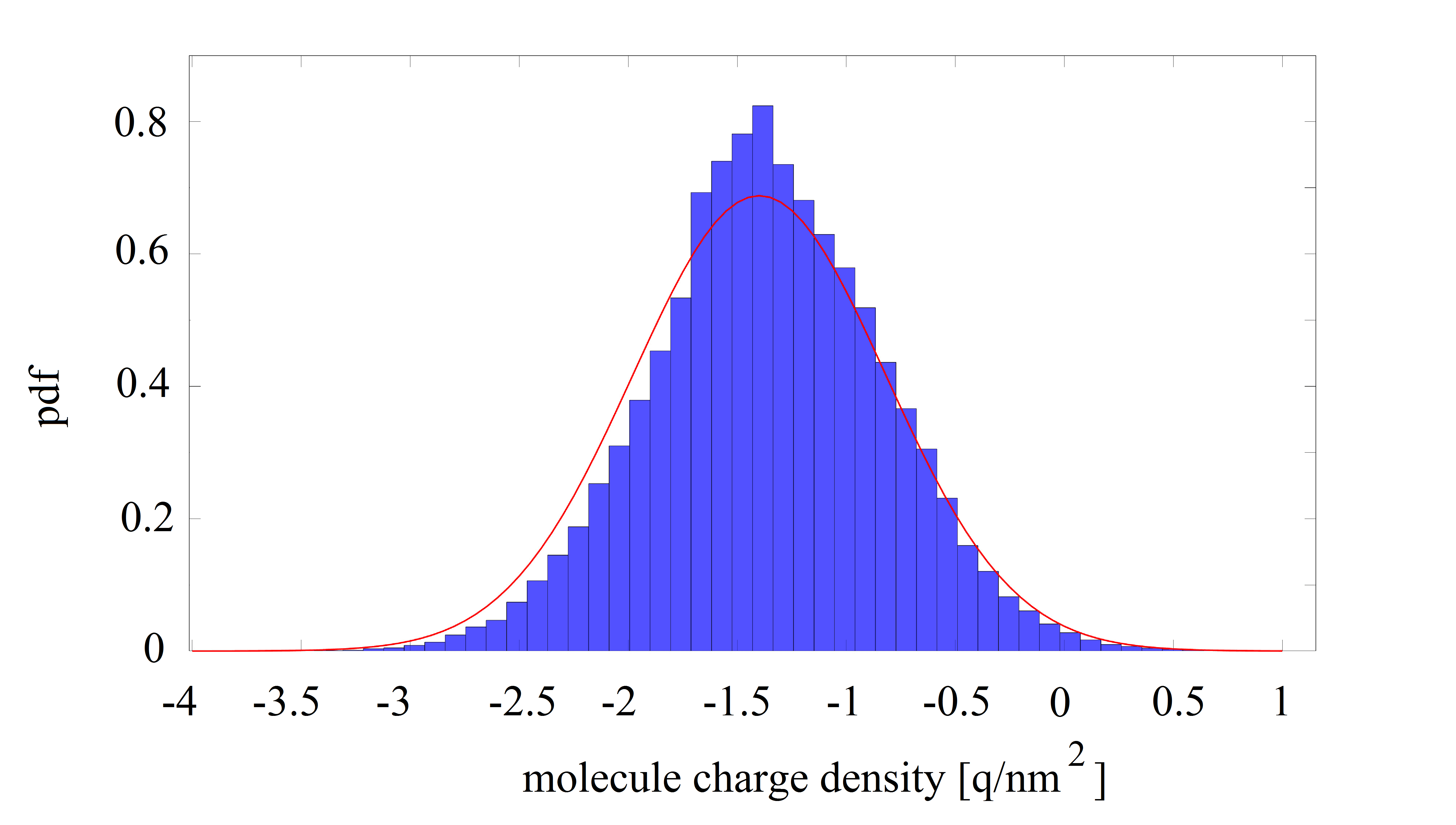}}
   \hfill    
   \subfloat{\includegraphics[width=8cm,height=5cm]{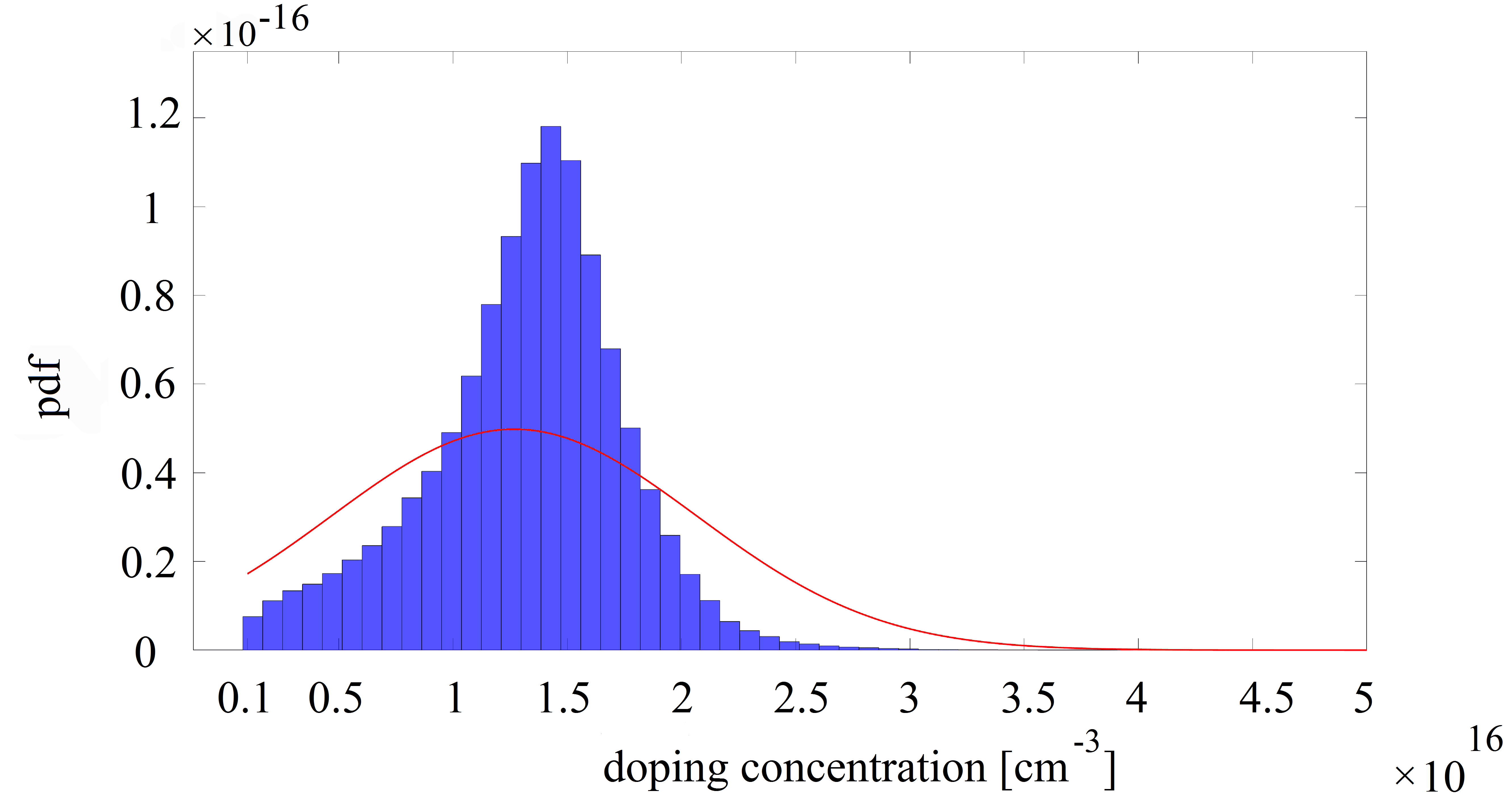}}
   \caption{The posterior (histogram) and prior (red line) distribution of the surface charge density (left) and the doping concentration (right) using the DRAM algorithm (marginal histograms).}
   \label{fig:hsit2}
 \end{figure} 

\begin{figure}[t!]
  \centering
  \subfloat{\includegraphics[width=8cm,height=5cm]{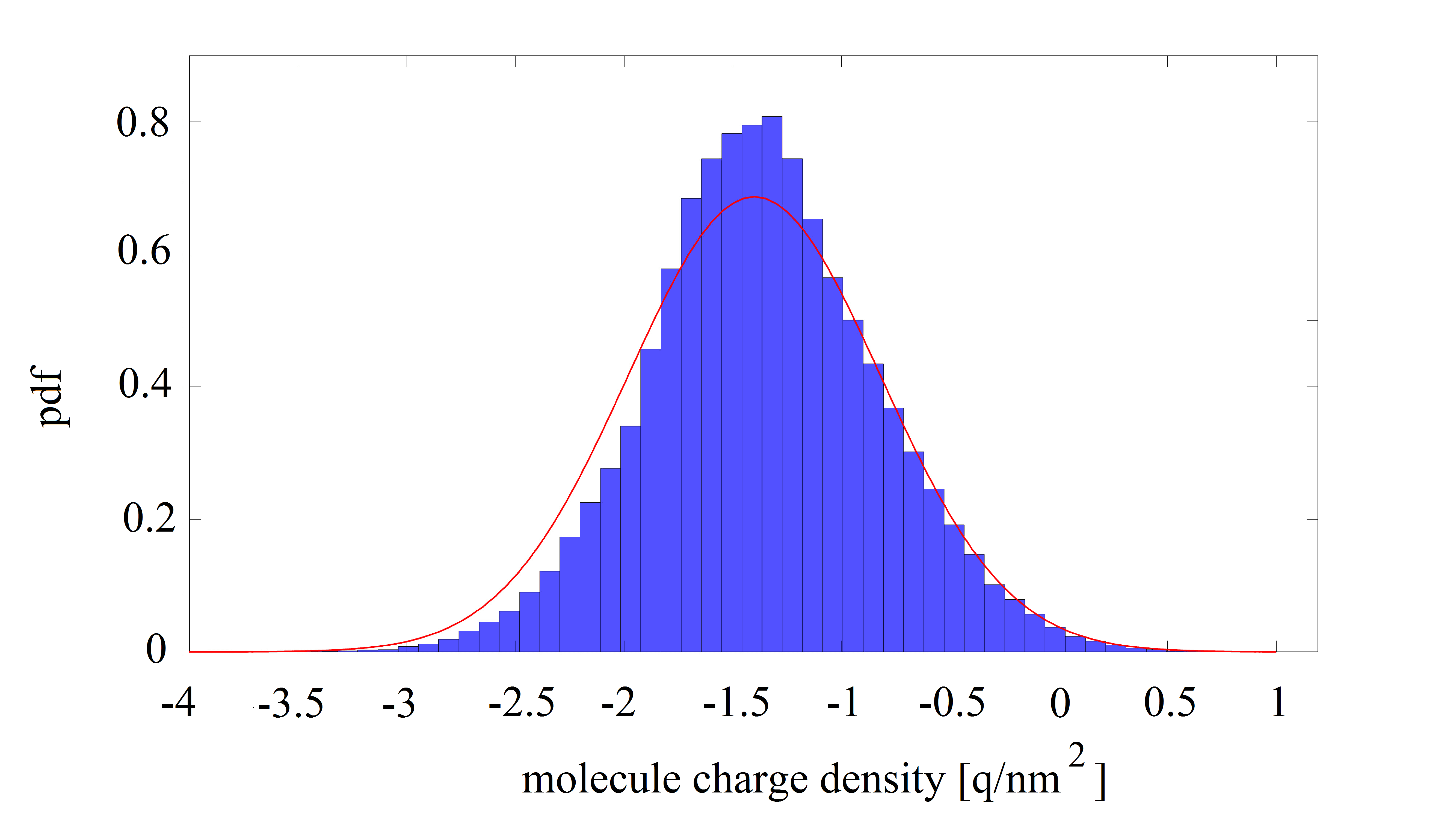}}%
  \hfill
  \subfloat{\includegraphics[width=8cm,height=5cm]{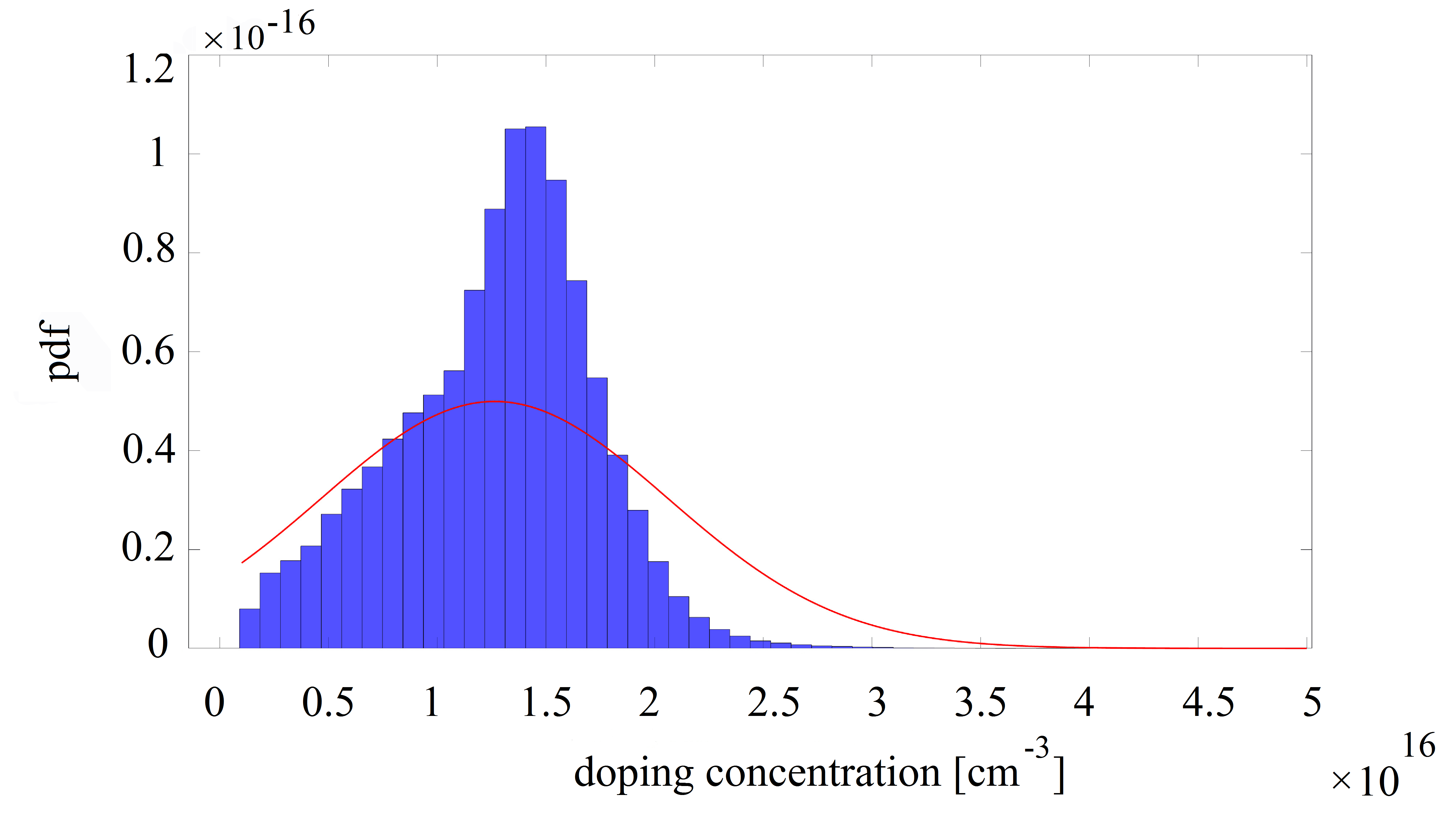}}%
  \newline
  \subfloat{\includegraphics[width=8cm,height=5cm]{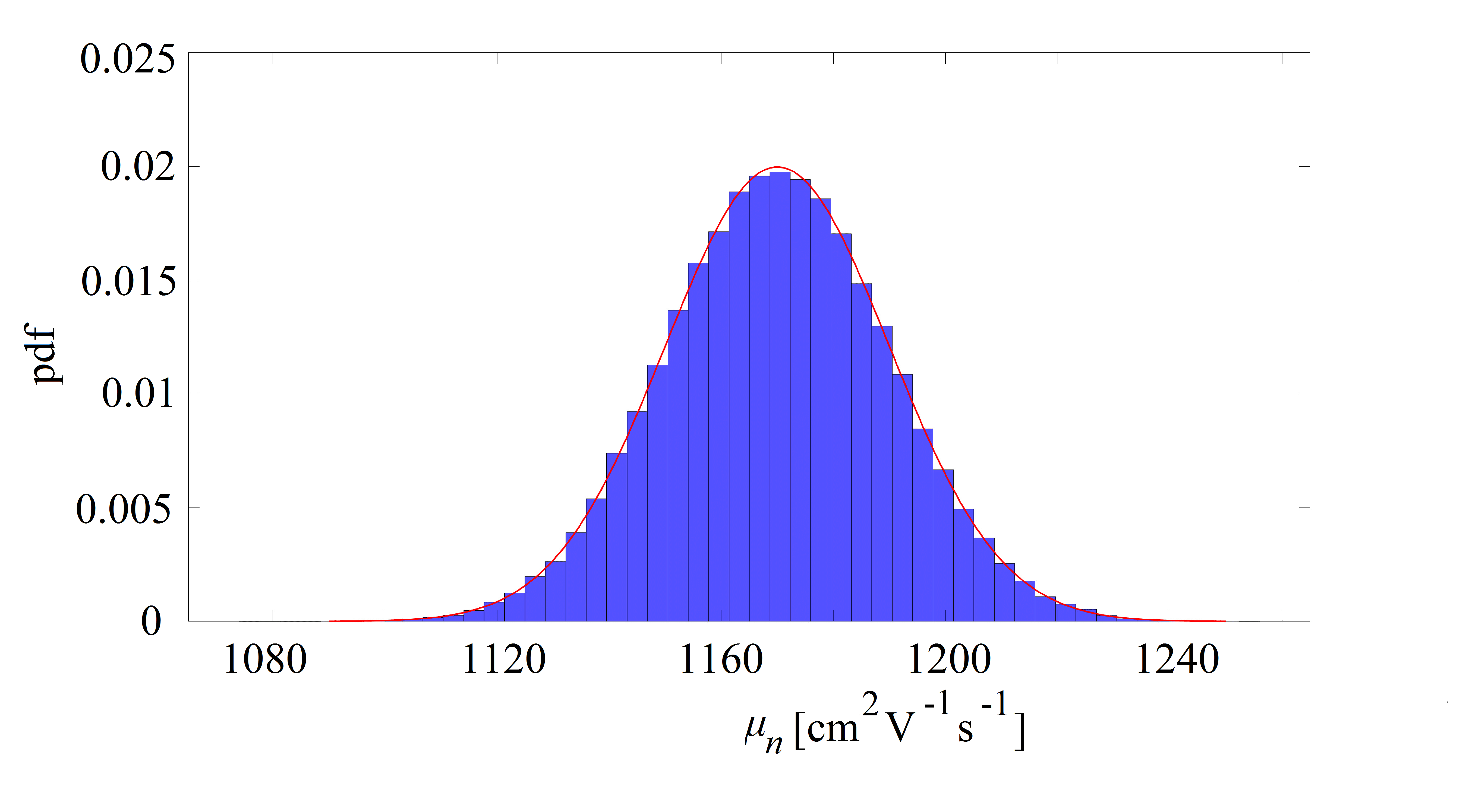}}   
  \hfill
    \subfloat{\includegraphics[width=8cm,height=5cm]{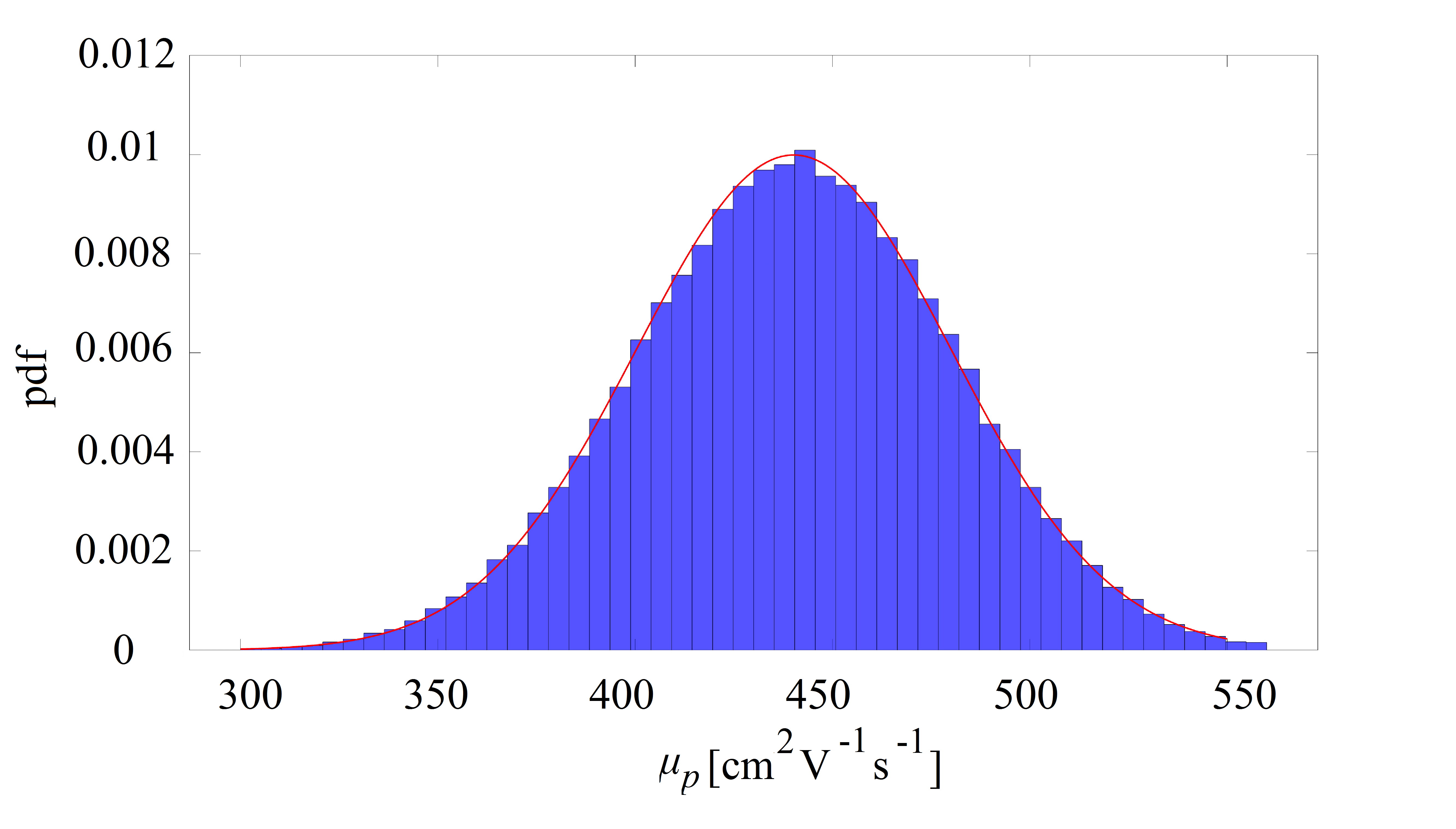}}
  \caption{The posterior distribution (histogram) and prior distribution (red line) of the surface charge density (top left), doping concentration (top right), electron mobility (bottom left) and  hole mobility (bottom right) using the DRAM algorithm.}
  \label{fig:hist3}
\end{figure}
  \begin{figure}[ht!]
 	\centering
 	\includegraphics[width=10cm,height=6cm]{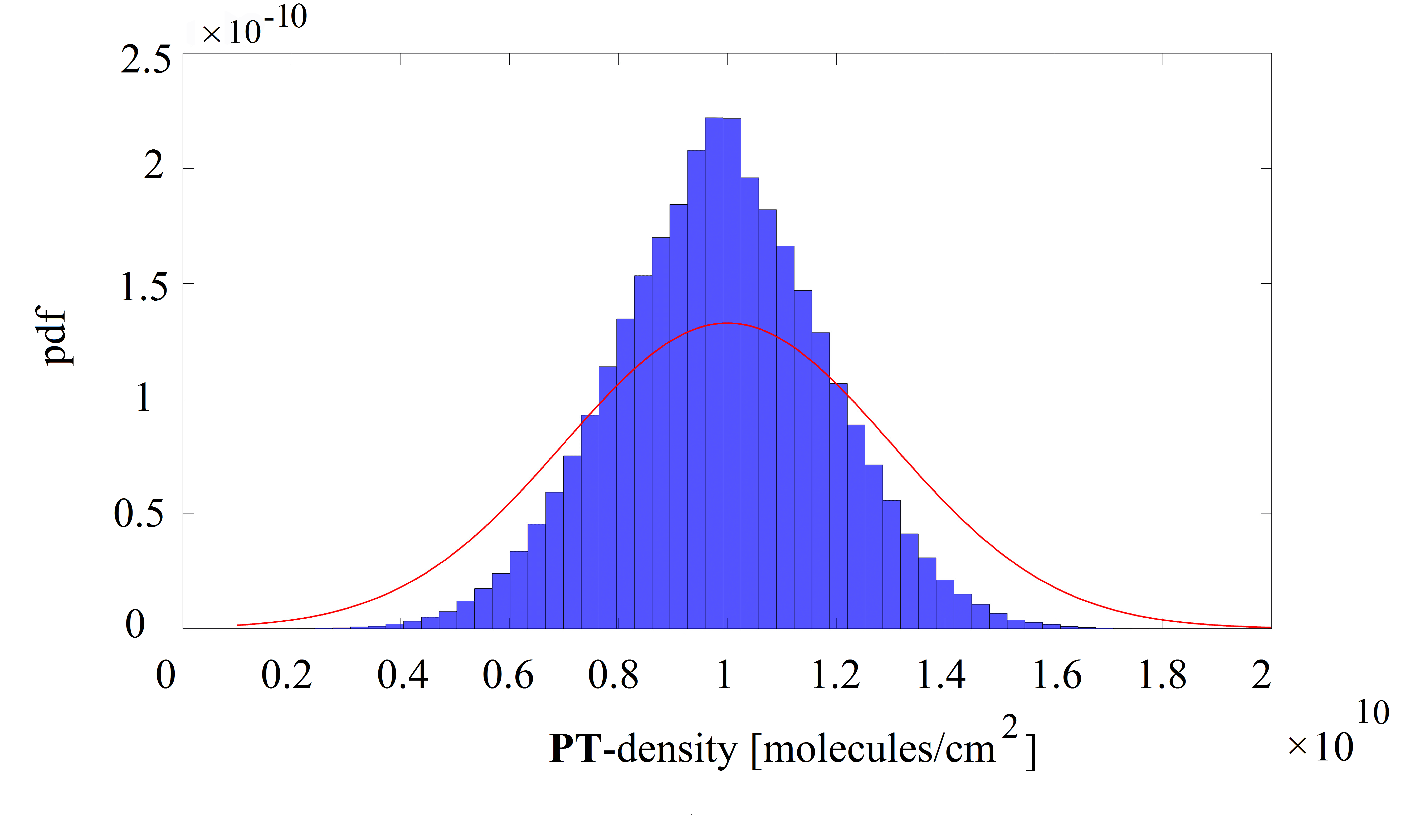}
 	\caption{The posterior distribution (histogram) and prior distribution (red line) of the \textbf{PT}-density using the DRAM algorithm.}
 	\label{fig:hist4}
 \end{figure}

 \begin{figure}[t!]
 	\centering
 	\subfloat{\includegraphics[width=8cm,height=5cm]{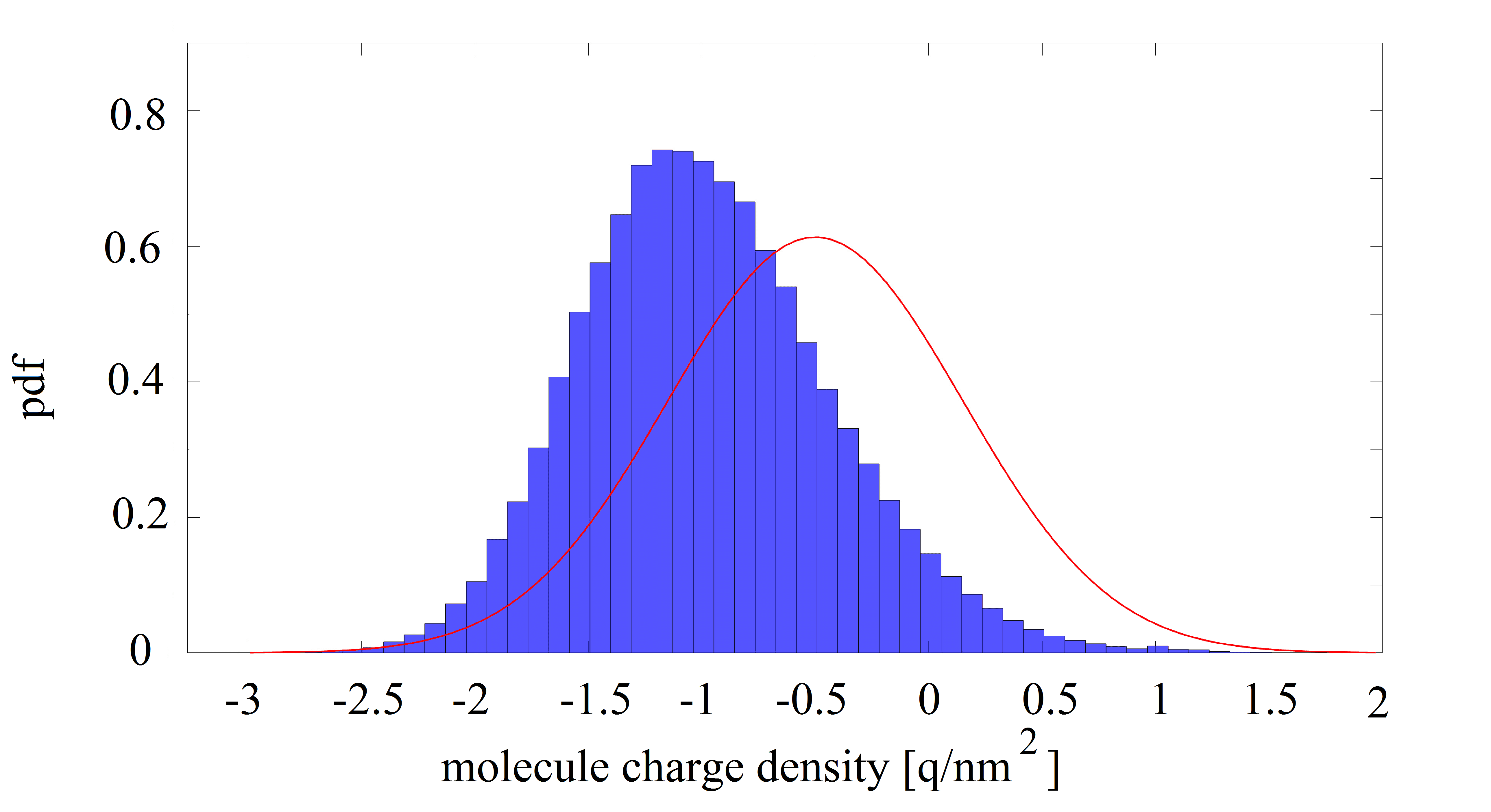}}%
 	\hfill
 	\subfloat{\includegraphics[width=8cm,height=5cm]{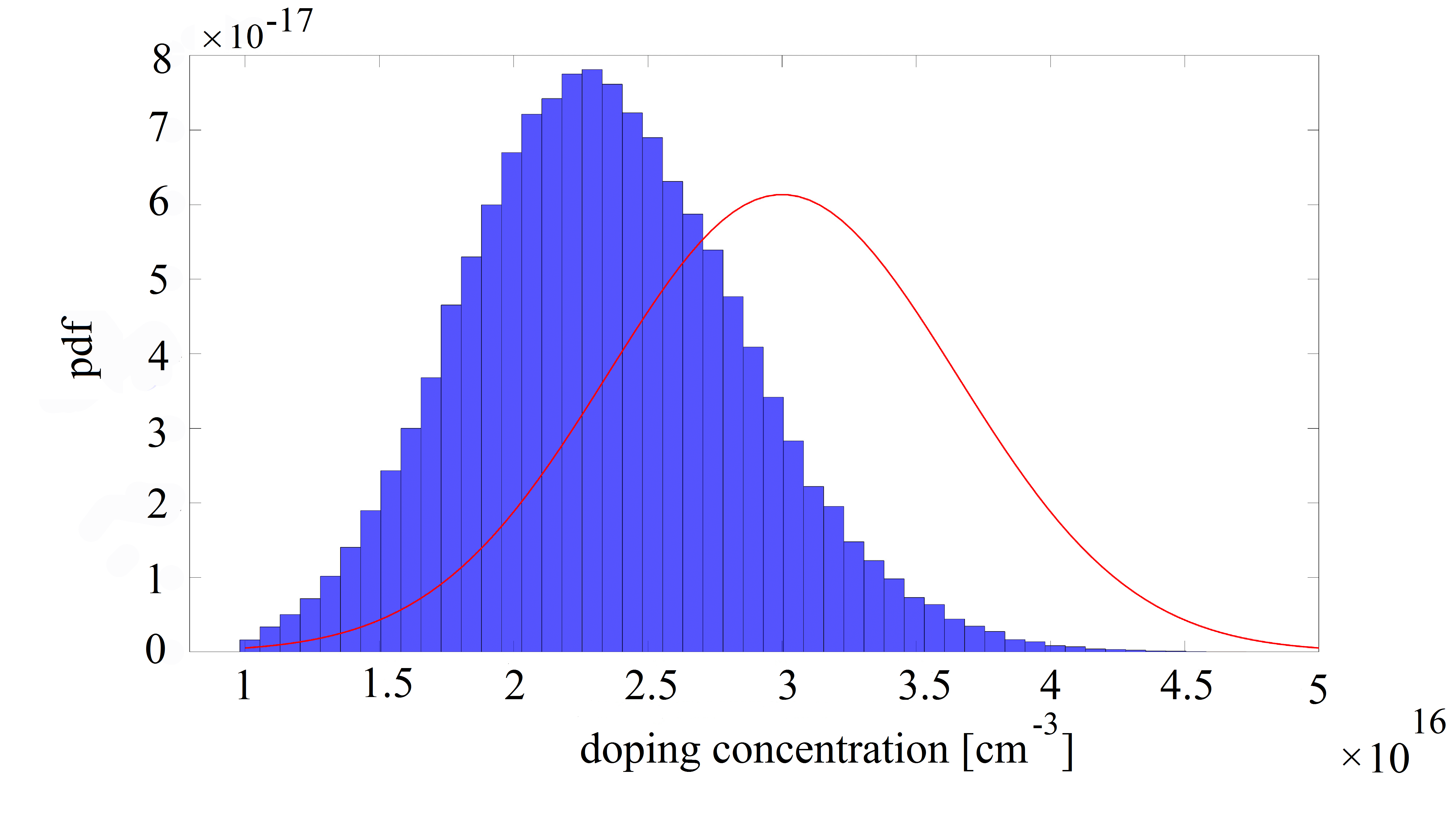}}%
 	\newline
 	\subfloat{\includegraphics[width=8cm,height=5cm]{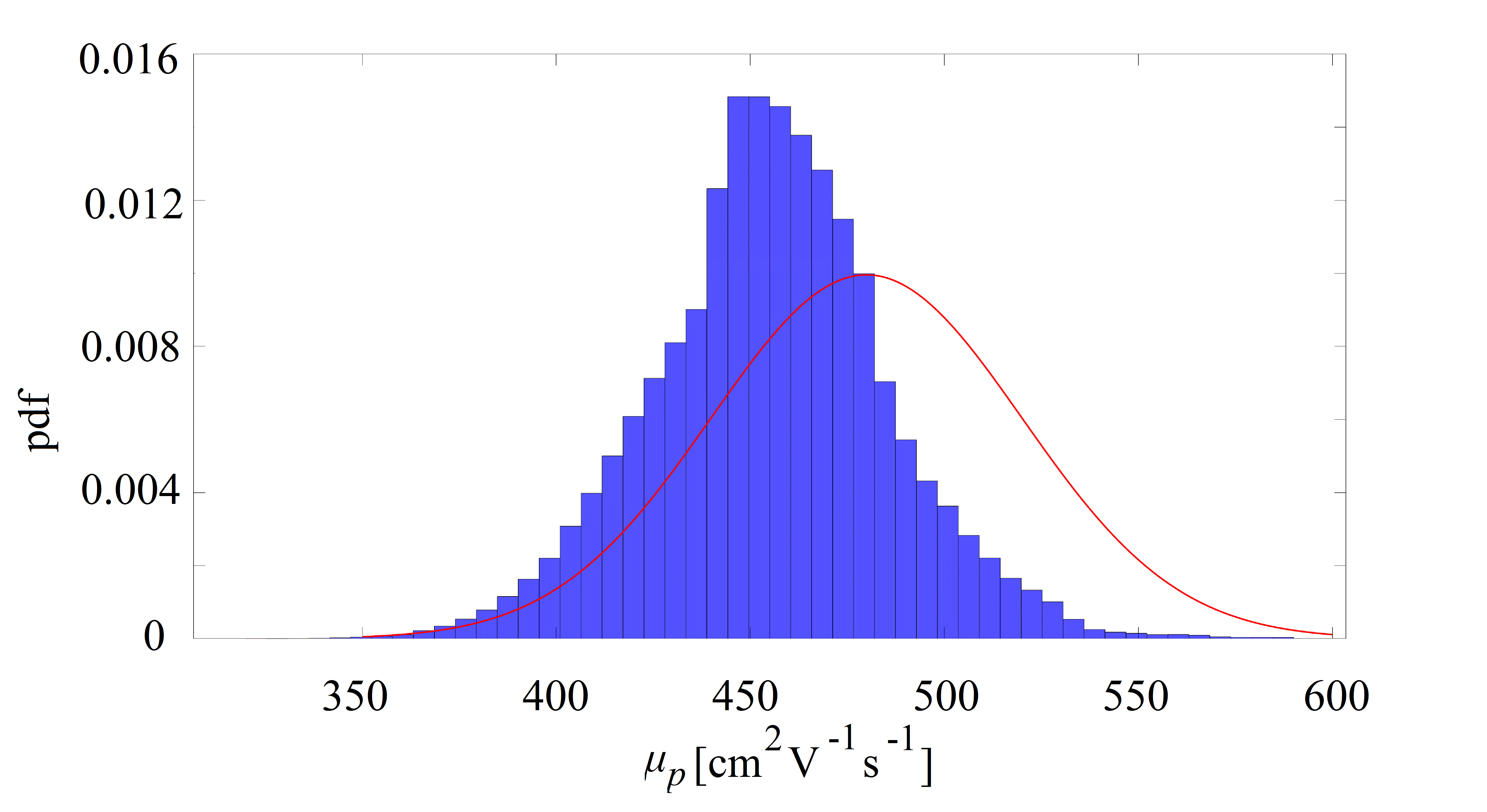}}   
 	\hfill
 	\subfloat{\includegraphics[width=8cm,height=5cm]{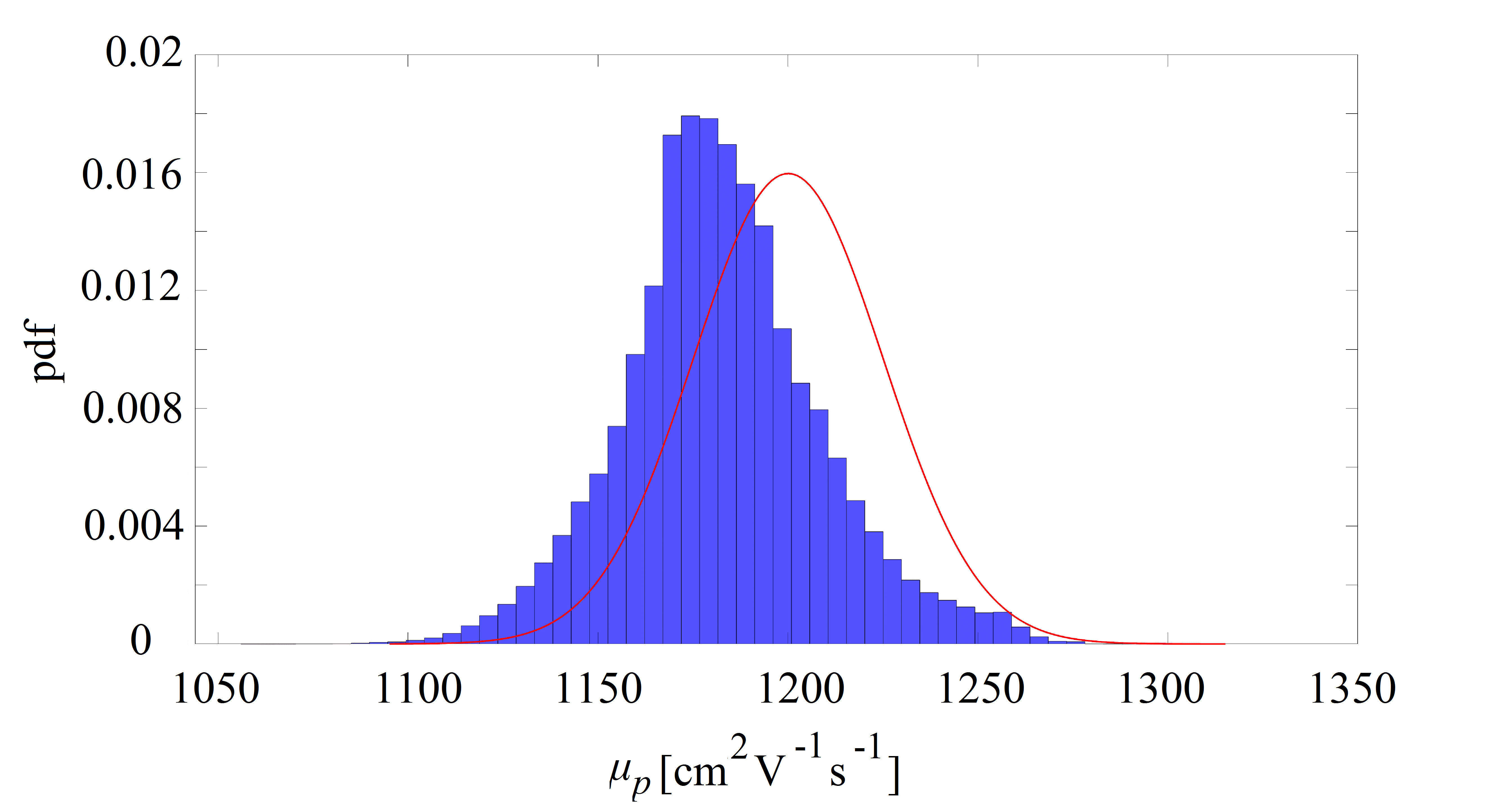}}
 	\caption{The posterior distribution (histogram) and prior distribution (red line) of the surface charge density (top left), doping concentration (top right), electron mobility (bottom left) and  hole mobility (bottom right) using the DRAM algorithm.}
 	\label{fig:hist}
 \end{figure}

In the next case, we simultaneously consider the effect of molecule charge density and doping concentration. In other words, the proposal $\alpha^*$ consists of two suggestions for the parameters (two-dimensional Bayesian estimation). In order to obtain the posterior distribution,  $\rho=\unit{-1.5}{q/nm^2}$ and $C_\text{dop}=\unit{10^{16}}{cm^{-3}}$ as the prior knowledge are applied in the simulations. We study the effect of molecule charge density from $\unit{1}{q/nm^2}$ to $\unit{-4}{q/nm^2}$ and of the doping concentration varying between $\unit{1\times10^{15}}{cm^{-3}}$ and $\unit{5\times10^{16}}{cm^{-3}}$. Figure \ref{fig:hsit2}  shows the posterior and prior (again Gaussian) distribution of doping concentration and molecule charge density where the acceptance rate is 62.3\%. Here, $\rho=\unit{-1.33}{q/nm^2}$ and $C_\text{dop}=\unit{1.18\times10^{16}}{cm^{-3}}$ as the expected values of the unknown parameters have been obtained. Similar to the first case, the probability of positive charge density is very low.


As we already mentioned, higher doping concentration increases the device conductivity; however, it decreases the sensor sensitivity. In other words, in high doping concentrations, $I$ (current with molecule) tends to $I_0$ \cite{khodadadian2017optimal} (current without molecules) since the effect of doping of the transducer is much more pronounced than the charged molecules.  
 Therefore, the chosen doping range strikes a balance between selectivity and conductivity.
As the figure shows, the most of accepted candidates are between $\unit{5\times10^{15}}{cm^{-3}}$ and  $\unit{2\times10^{15}}{cm^{-6}}$ and
the probability of $C_\text{dop}>\unit{3\times10^{16}}{cm^{-3}}$ is negligible, which confirms the doping effect on the sensitivity.


The electron and hole mobilities have a similar dependence on doping. For low doping concentrations, the mobility
is almost constant and primarily limited by phonon scattering. At higher doping concentrations the mobility
decreases due to ionized impurity scattering with the ionized doping atoms. In \cite{arora1982electron}, an analytic expression for electron and hole mobility in silicon as a function of doping concentration has been given. In the previous cases we used the electron/hole mobility according to the Arora formula, i.e., $\mu_p=\unit{430}{cm^{-2}V^{-1}s^-1}$ and $\mu_n=\unit{1170}{cm^{-2}V^{-1}s^-1}$ for $C_\text{dop}=\unit{1\times10^{16}}{cm^{-3}}$. Now in order to validate this empirical formula, we consider the mobilities as the other unknown parameters. Figure \ref{fig:hist3} illustrates the posterior distribution of four physical parameters, where $\rho=\unit{-1.38}{q/nm^2}$, $C_\text{dop}=\unit{1.91\times10^{16}}{cm^{-3}}$,  $\mu_n=\unit{1175}{cm^{-2}V^{-1}s^-1}$ and $\mu_p=\unit{439}{cm^{-2}V^{-1}s^-1}$  are found as the expected values. The obtained mobilities also confirm the Arora formula. Again, the (Gaussian) prior distribution is shown, and for this estimation, the acceptance rate of 59.8 \% achieved.

We consider probe-target binding in the equilibrium
and use a receptor concentration of $C_P=\unit{3\times 10^{11}}{cm^{-2}}$. In practice, a good estimation of the number of bound target to the probe molecules cannot be achieved easily. Here we study the density between $1\times 10^{9}$ and $2\times 10^{10}$ molecules per square centimeters (the Gaussian prior distribution). 
Figure \ref{fig:hist4} shows the density estimation where the rest of (four) unknown parameters are according to the extracted information by the posterior distributions.  As shown, the mean of $\bf{PT}$-density is $\unit{1.05\times10^{10}}{mol/cm^{2}}$ and 75.8\% is the acceptance rate.

As we already mentioned, using reliable prior knowledge (good guess) by employing PROPKA algorithm, Arora formula and good approximation of doping density enables us to provide an efficient posterior distribution. This fact gives rise to uncertainty reduction of parameters and also a good acceptance rate is achieved. Now we study the effect of the prior distribution on the marginal posterior. To that end, for molecule charge density, we assume it varies between $\unit{-3}{q/nm^2}$ and $\unit{2}{q/nm^2}$ (with the mean of $\unit{-0.5}{q/nm^2}$), the doping concentration changes from $\unit{1\times 10^{16}}{cm^{-3}}$ to $\unit{5\times 10^{16}}{cm^{-3}}$ (the mean is $\unit{2.5\times 10^{16}}{cm^{-3}}$) and regarding the mobilities, $\mu_p=\unit{480}{cm^{-2}V^{-1}s^{-1}}$ and $\mu_n=\unit{1200}{cm^{-2}V^{-1}s^{-1}}$ are chosen. Using the mentioned proposals lead to the mean value of $\rho=\unit{-0.9577}{q/nm^2}$,  $C_\text{dop}=\unit{2.235\times 10^{16}}{cm^{-3}}$,  $\mu_p=\unit{455}{cm^{-2}V^{-1}s^-1}$ and $\mu_n=\unit{1190}{cm^{-2}V^{-1}s^-1}$. As a noticeable difference with the previous estimation, in spite of the convergence, the acceptance rate reduced to 32.3 \%.

\begin{figure}[ht!]
	\centering
	\includegraphics[width=10cm,height=7cm]{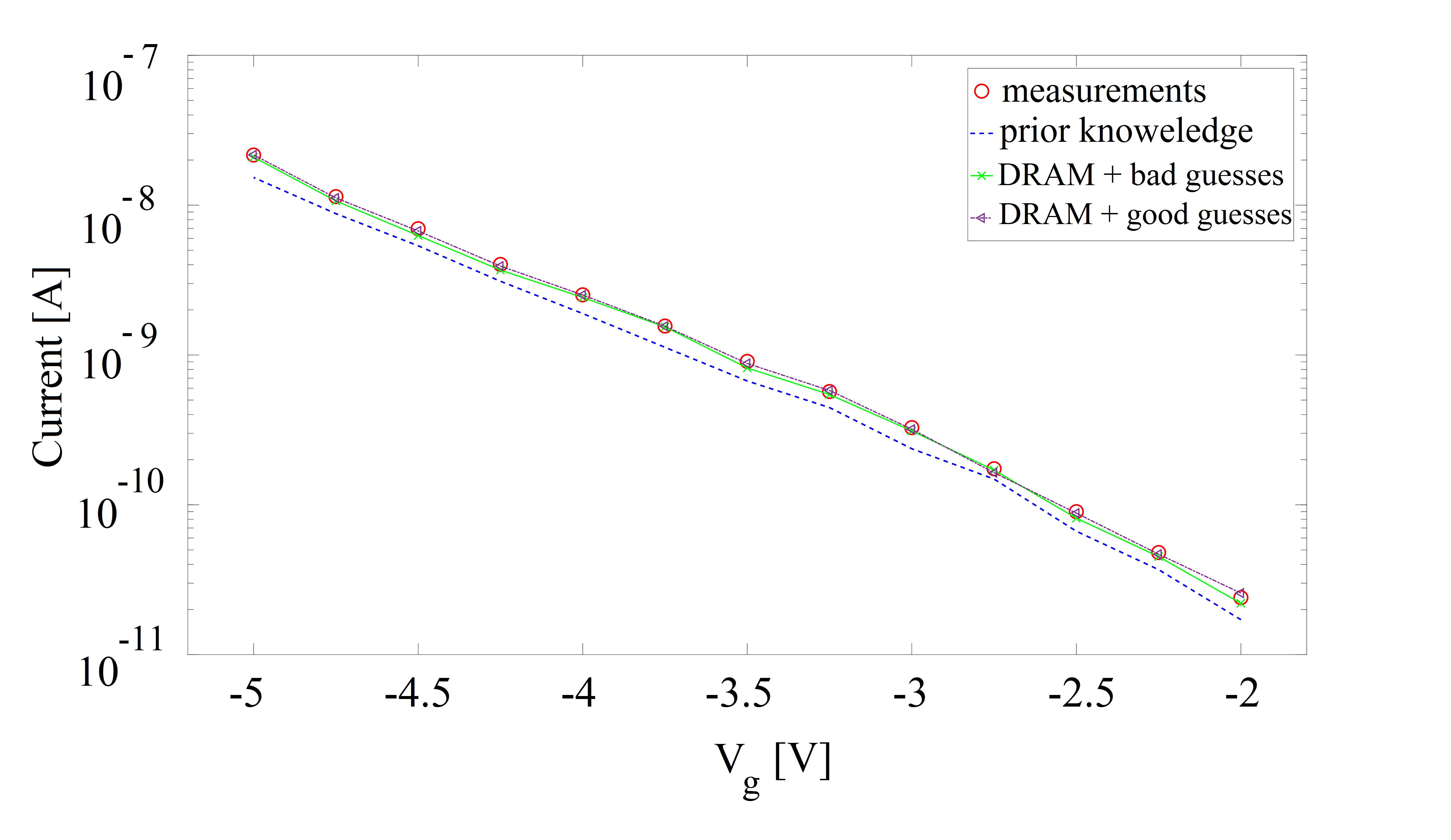}
	\caption{The experimental data versus simulated currents. Here, we calculated the current according to the prior information and the extracted posterior knowledge. The I-V curve shows the significant advantage of the Bayesian inference where the simulations have a very good agreement with the measurements.}
	\label{fig:IV2}
\end{figure}
Now we employ the \blue{parameters estimated by taking the mean values
  of their posterior distributions} to calculate the electrical current. We have obtained two posterior distributions first based on the empirical formulas and second according to not good guesses.  
Figure \ref{fig:IV2} shows the simulated current as a function of different gate voltages for both posterior distributions and compared it with the experiments.  These results validate the effectiveness and usefulness of the Bayesian inference since using the DRAM algorithm leads to an excellent agreement between the measurement and the simulation. However, the posterior distribution with a reasonable guess gives rise to a more exact electrical current.


\begin{figure}[ht!]
   \centering
   \subfloat{\includegraphics[width=8cm,height=5cm]{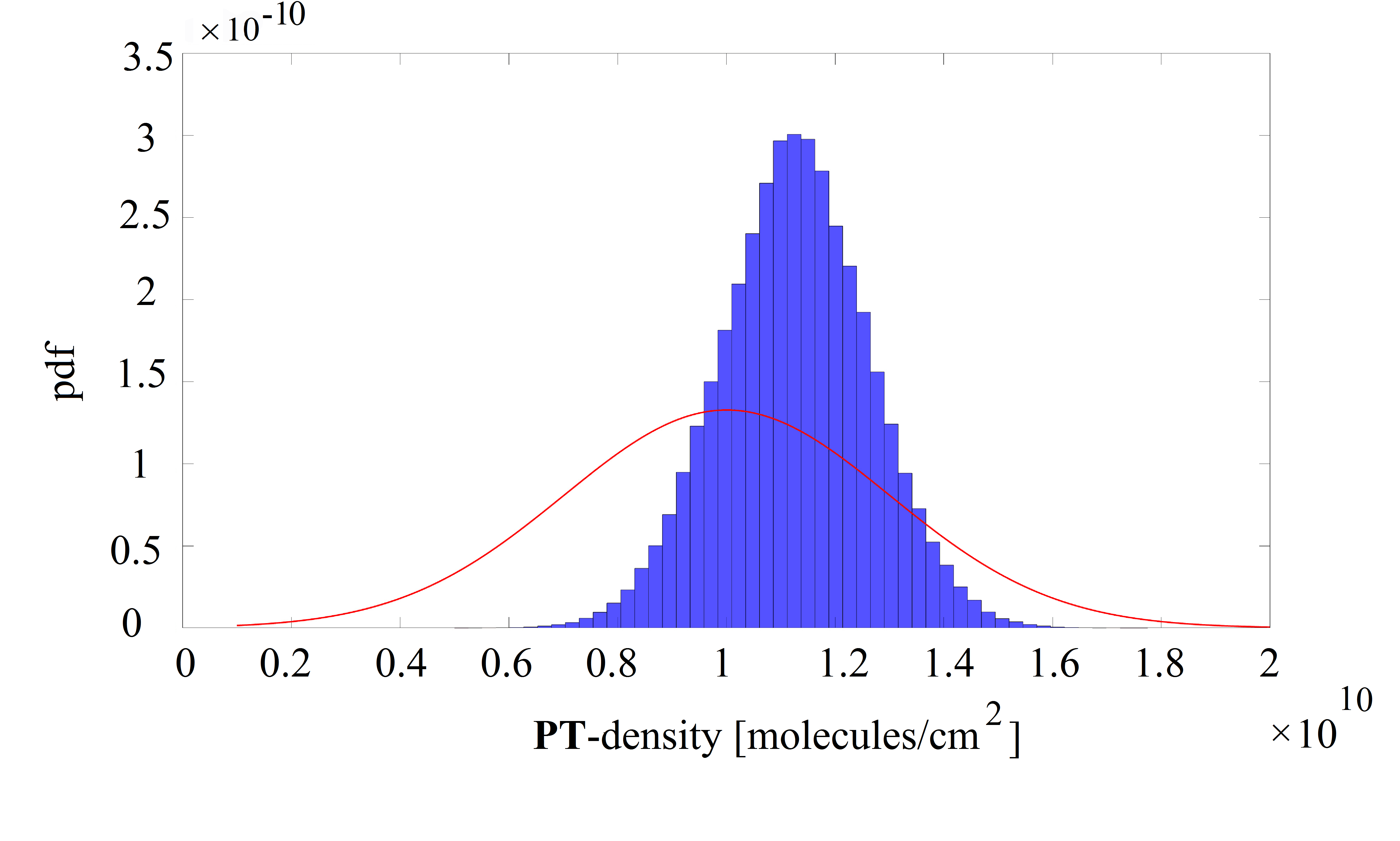}}
   \hfill    
   \subfloat{\includegraphics[width=8cm,height=5cm]{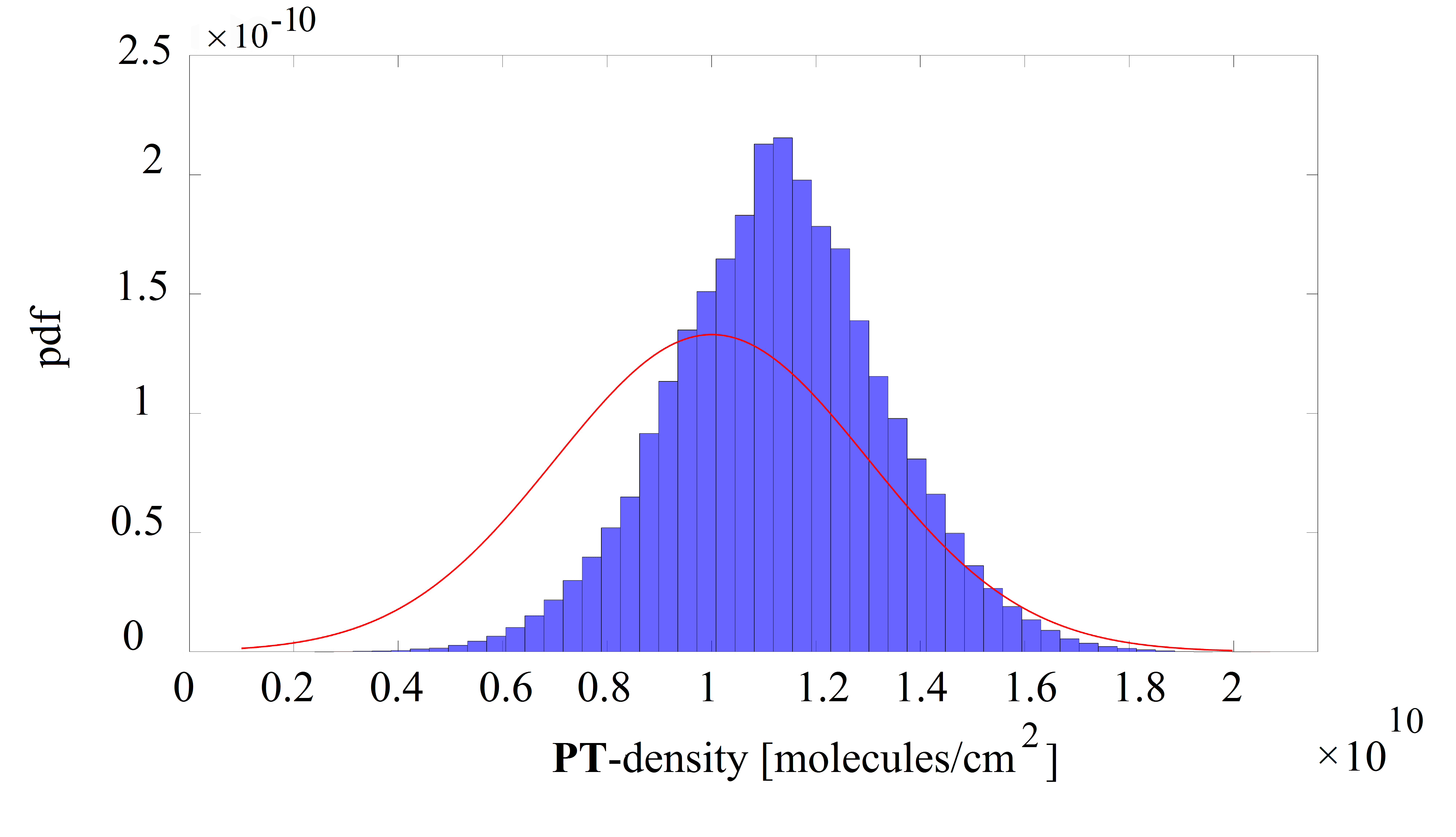}}
   \caption{The posterior (histogram) and prior (red line) distribution  of the $\bf{PT}-$complex density using $\rho=\unit{-1.38}{q/nm^2}$, $C_\text{dop}=\unit{1.91\times10^{16}}{cm^{-3}}$, $\mu_n=\unit{1175}{cm^{-2}V^{-1}s^{-1}}$, and   $\mu_p=\unit{439}{cm^{-2}V^{-1}s^{-1}}$ with $\varepsilon=0.5\%$ (left) and $\varepsilon=2\%$ (right).}
   \label{fig:PT1}
 \end{figure}

The measurement error $\varepsilon$ affects the Bayesian inversion as well. We note that in addition to the measurement error, the spatial discretization error  (finite element discretization) and the statistical error (Monte Carlo sampling) might be part of $\varepsilon$. Here, we study the influence of $\varepsilon$ on the $\bf{PT}$-complex density. Figure \ref{fig:PT1} shows the posterior distribution for two more measurement errors, namely $\varepsilon=0.5\%$ and $\varepsilon=2\%$. Also, the distribution for $\varepsilon=1\%$ was already illustrated in Figure \ref{fig:hist4}. The results point out that a smaller error gives rise to a narrower distribution, while a larger error makes it wider, as expected.

\section{Conclusions}


 In  sensor design, reliable information about different device parameters is crucial. The 
stochastic DDPB system is a useful computational system to model the electrical/electrochemical behavior of nanowire sensors. The model enables us to study the effect of different influential parameters, e.g., molecule charge density, diffusion coefficients, doping concentration, gate voltage, etc.
We have used Bayesian inversion to provide a reliable estimation of the parameters. More precisely, the MCMC method (DRAM algorithm) has been used to obtain the  posterior distribution of the physical device parameters.

In this paper, we have first validated the DDPB system with the experimental data and then applied the DRAM algorithm to estimate the parameters. In order to study the effect of charged molecules on the nanowire, we estimated the molecule charge density. Then, we considered the effect of doping concentration (in a two-dimensional Bayesian inversion). Here, due to the narrow probability density for each parameter,  reliable information can be extracted. 
  In addition to the mentioned physical parameters, we studied the effect of electron and hole mobilities in addition to the previous unknowns (in a four-dimensional Bayesian estimation) and provided their posterior distributions simultaneously. In the most complicated simulation, we have estimated the probability density of the $\bf{PT}$-concentration. The results enable us to determine the device and molecule properties at the same time. 

Finally, we have applied the results obtained by Bayesian inference to the DDPB system and again simulated the device current. The results point out that compared to the previous simulations, the agreement with the experimental data has improved, which indicates the effectiveness of the DRAM technique.  The results show that Bayesian inversion is a promising technique and
has  significant capabilities in the design of various sensors and nanoscale devices as well as in interpreting measurement data and assessing its quality.

\label{conclusions}

\section{Acknowledgments}     

The  authors acknowledge support by the FWF (Austrian Science Fund) START
project No.\ Y660 \textit{PDE Models for Nanotechnology}. \blue{The authors also acknowledge the helpful
comments by the anonymous reviewers.}
  
%

 


\begin{thebibliography}{10}
	\expandafter\ifx\csname url\endcsname\relax
	\def\url#1{\texttt{#1}}\fi
	\expandafter\ifx\csname urlprefix\endcsname\relax\def\urlprefix{URL }\fi
	\expandafter\ifx\csname href\endcsname\relax
	\def\href#1#2{#2} \def\path#1{#1}\fi
	
	\bibitem{duan2013complementary}
	X.~Duan, N.~K. Rajan, M.~H. Izadi, M.~A. Reed, Complementary metal oxide
	semiconductor-compatible silicon nanowire biofield-effect transistors as
	affinity biosensors, Nanomedicine 8~(11) (2013) 1839--1851.
	
	\bibitem{cui2001nanowire}
	Y.~Cui, Q.~Wei, H.~Park, C.~M. Lieber, Nanowire nanosensors for highly
	sensitive and selective detection of biological and chemical species, Science
	293~(5533) (2001) 1289--1292.
	
	\bibitem{zheng2005multiplexed}
	G.~Zheng, F.~Patolsky, Y.~Cui, W.~U. Wang, C.~M. Lieber, Multiplexed electrical
	detection of cancer markers with nanowire sensor arrays, Nature Biotechnology
	23~(10) (2005) 1294.
	
	\bibitem{he2015label}
	J.~He, J.~Zhu, C.~Gong, J.~Qi, H.~Xiao, B.~Jiang, Y.~Zhao, Label-free direct
	detection of {miRNAs with} poly-silicon nanowire biosensors, PLoS One 10~(12)
	(2015) e0145160.
	
	\bibitem{hahm2004direct}
	J.-i. Hahm, C.~M. Lieber, Direct ultrasensitive electrical detection of{ DNA
		and DNA} sequence variations using nanowire nanosensors, Nano Letters 4~(1)
	(2004) 51--54.
	
	\bibitem{wang2005label}
	W.~U. Wang, C.~Chen, K.-h. Lin, Y.~Fang, C.~M. Lieber, Label-free detection of
	small-molecule--protein interactions by using nanowire nanosensors,
	Proceedings of the National Academy of Sciences 102~(9) (2005) 3208--3212.
	
	\bibitem{chua2009label}
	J.~H. Chua, R.-E. Chee, A.~Agarwal, S.~M. Wong, G.-J. Zhang, Label-free
	electrical detection of cardiac biomarker with complementary metal-oxide
	semiconductor-compatible silicon nanowire sensor arrays, Analytical Chemistry
	81~(15) (2009) 6266--6271.
	
	\bibitem{baumgartner2013predictive}
	S.~Baumgartner, C.~Heitzinger, A.~Vacic, M.~A. Reed, {Predictive simulations
		and optimization of nanowire field-effect PSA sensors including screening},
	Nanotechnology 24~(22) (2013) 225503.
	
	\bibitem{shashaani2016silicon}
	H.~Shashaani, M.~Faramarzpour, M.~Hassanpour, N.~Namdar, A.~Alikhani,
	M.~Abdolahad, Silicon nanowire based biosensing platform for electrochemical
	sensing of {Mebendazole} drug activity on breast cancer cells, Biosensors and
	Bioelectronics 85 (2016) 363--370.
	
	\bibitem{lee2010measurements}
	M.-H. Lee, D.-H. Lee, S.-W. Jung, K.-N. Lee, Y.~S. Park, W.-K. Seong,
	{Measurements of serum C-reactive protein levels in patients with gastric
		cancer and quantification using silicon nanowire arrays}, Nanomedicine:
	Nanotechnology, Biology and Medicine 6~(1) (2010) 78--83.
	
	\bibitem{shen2012rapid}
	F.~Shen, J.~Wang, Z.~Xu, Y.~Wu, Q.~Chen, X.~Li, X.~Jie, L.~Li, M.~Yao, X.~Guo,
	et~al., Rapid flu diagnosis using silicon nanowire sensor, Nano Letters
	12~(7) (2012) 3722--3730.
	
	\bibitem{guan2014highly}
	W.~Guan, X.~Duan, M.~A. Reed, Highly specific and sensitive non-enzymatic
	determination of uric acid in serum and urine by extended gate field effect
	transistor sensors, Biosensors and Bioelectronics 51 (2014) 225--231.
	
	\bibitem{mu2015silicon}
	L.~Mu, Y.~Chang, S.~D. Sawtelle, M.~Wipf, X.~Duan, M.~A. Reed, {Silicon
		nanowire field-effect transistors—a versatile class of potentiometric
		nanobiosensors}, IEEE Access 3 (2015) 287--302.
	
	\bibitem{khodadadian2017optimal}
	A.~Khodadadian, K.~Hosseini, A.~Manzour-ol Ajdad, M.~Hedayati,
	R.~Kalantarinejad, C.~Heitzinger, Optimal design of nanowire field-effect
	troponin sensors, Computers in Biology and Medicine 87 (2017) 46--56.
	
	\bibitem{mirsian2019new}
	S.~Mirsian, A.~Khodadadian, M.~Hedayati, A.~Manzour-ol Ajdad,
	R.~Kalantarinejad, C.~Heitzinger, A new method for selective
	functionalization of silicon nanowire sensors and {Bayesian} inversion for
	its parameters, Biosensors and Bioelectronics 142 (2019) 111527.
	
	\bibitem{khodadadian2016basis}
	A.~Khodadadian, C.~Heitzinger, {Basis adaptation for the stochastic nonlinear
		Poisson--Boltzmann equation}, Journal of Computational Electronics 15~(4)
	(2016) 1393--1406.
	
	\bibitem{taghizadeh2017optimal}
	L.~Taghizadeh, A.~Khodadadian, C.~Heitzinger, {The optimal multilevel
		Monte-Carlo approximation of the stochastic drift--diffusion-Poisson system},
	Computer Methods in Applied Mechanics and Engineering 318 (2017) 739--761.
	
	\bibitem{khodadadian2018three}
	A.~Khodadadian, L.~Taghizadeh, C.~Heitzinger, {Three-dimensional optimal
		multi-level Monte--Carlo approximation of the stochastic
		drift--diffusion--Poisson system in nanoscale devices}, Journal of
	Computational Electronics 17~(1) (2018) 76--89.
	
	\bibitem{khodadadian2018optimal}
	A.~Khodadadian, L.~Taghizadeh, C.~Heitzinger, {Optimal multilevel randomized
		quasi-Monte-Carlo method for the stochastic drift--diffusion-Poisson system},
	Computer Methods in Applied Mechanics and Engineering 329 (2018) 480--497.
	
	\bibitem{dashti2017bayesian}
	M.~Dashti, A.~M. Stuart, {The Bayesian approach to inverse problems}, Handbook
	of Uncertainty Quantification (2017) 311--428.
	
	\bibitem{smith2013uncertainty}
	R.~C. Smith, Uncertainty Quantification: Theory, Implementation, and
	Applications, Vol.~12, SIAM, 2013.
	
	\bibitem{smith1993bayesian}
	A.~F. Smith, G.~O. Roberts, {Bayesian computation via the Gibbs sampler and
		related Markov chain Monte Carlo methods}, Journal of the Royal Statistical
	Society. Series B (Methodological) (1993) 3--23.
	
	\bibitem{andrieu2003introduction}
	C.~Andrieu, N.~De~Freitas, A.~Doucet, M.~I. Jordan, {An introduction to MCMC
		for machine learning}, Machine Learning 50~(1-2) (2003) 5--43.
	
	\bibitem{mira2001metropolis}
	A.~Mira, et~al., {On Metropolis-Hastings algorithms with delayed rejection},
	Metron 59~(3-4) (2001) 231--241.
	
	\bibitem{haario2006dram}
	H.~Haario, M.~Laine, A.~Mira, E.~Saksman, {DRAM: efficient adaptive MCMC},
	Statistics and Computing 16~(4) (2006) 339--354.
	
	\bibitem{zuev2011modified}
	K.~M. Zuev, L.~S. Katafygiotis, {Modified Metropolis--Hastings algorithm with
		delayed rejection}, Probabilistic Engineering Mechanics 26~(3) (2011)
	405--412.
	
	\bibitem{punzet2012determination}
	M.~Punzet, D.~Baurecht, F.~Varga, H.~Karlic, C.~Heitzinger, {Determination of
		surface concentrations of individual molecule-layers used in nanoscale
		biosensors by in situ ATR-FTIR spectroscopy}, Nanoscale 4~(7) (2012)
	2431--2438.
	
	\bibitem{heitzinger2010multiscale}
	C.~Heitzinger, N.~J. Mauser, C.~Ringhofer, Multiscale modeling of planar and
	nanowire field-effect biosensors, SIAM Journal on Applied Mathematics 70~(5)
	(2010) 1634--1654.
	
	\bibitem{tulzer2014fluctuations}
	G.~Tulzer, C.~Heitzinger, Fluctuations due to association and dissociation
	processes at nanowire-biosensor surfaces and their optimal design,
	Nanotechnology 26~(2) (2014) 025502.
	
	\bibitem{stuart2010inverse}
	A.~M. Stuart, {Inverse problems: a Bayesian perspective}, Acta Numerica 19
	(2010) 451--559.
	
	\bibitem{li2005very}
	H.~Li, A.~D. Robertson, J.~H. Jensen, {Very fast empirical prediction and
		rationalization of protein pKa values}, Proteins: Structure, Function, and
	Bioinformatics 61~(4) (2005) 704--721.
	
	\bibitem{arora1982electron}
	N.~D. Arora, J.~R. Hauser, D.~J. Roulston, Electron and hole mobilities in
	silicon as a function of concentration and temperature, IEEE Transactions on
	electron devices 29~(2) (1982) 292--295.
	
\end{thebibliography}


\end{document}